\documentclass[11pt]{article}

\usepackage{latexsym,graphicx,epstopdf,amsmath,varioref,algorithm,algorithmicx,algpseudocode,amsfonts}
\DeclareGraphicsExtensions{.eps,.jpg}

\newcommand{\beqn}[1]{\begin{equation}\label{#1}}
\newcommand{\eeqn}{\end{equation}}

\newcommand{\eqdef}{\stackrel{\rm def}{=}}

\newcommand{\bctable}[1]{\begin{table}[htbp]
                         \begin{center}
                         \begin{tabular}{#1} }
\newcommand{\ectable}[1]{\end{tabular}
                         \caption{#1}
                         \end{center}
                         \end{table}}
\newtheorem{theorem}{Theorem}[section]
\newtheorem{lemma}[theorem]{Lemma}

\newcommand{\hatq}{\hat{q}}

\newcommand{\comment}[1]{}

\newcommand{\normA}[1]{\| #1 \|_A}
\newcommand{\normAi}[1]{\| #1 \|_{A^{-1}}}
\newcommand{\normAiA}[1]{\| #1 \|_{A^{-1},A}}

\usepackage{color}
\usepackage{color}
\usepackage{mathtools}

\title{Minimizing convex quadratics with variable precision conjugate gradients}

\author{
  S. Gratton
  \thanks{Universit\'e de Toulouse, INP, IRIT, Toulouse, France, serge.gratton@toulouse-inp.fr},
  E. Simon
  \thanks{Universit\'e de Toulouse, INP, IRIT, Toulouse, France, ehouarn.simon@toulouse-inp.fr},
  D. Titley-Peloquin
	\thanks{McGill University, Montreal, Canada, david.titley-peloquin@mcgill.ca},
  and Ph. L. Toint
  \thanks{NAXYS, University of Namur, Namur, Belgium, philippe.toint@unamur.be.
          Partially supported by ANR-11-LABX-0040-CIMI within the program ANR-11-IDEX-0002-02.}
}

\begin{document}

\maketitle

\begin{abstract}
We investigate the method of conjugate gradients,
exploiting inaccurate matrix-vector products,
for the solution of convex quadratic optimization problems.
Theoretical performance bounds are derived,
and the necessary quantities occurring in the theoretical bounds estimated,
leading to a practical algorithm.
Numerical experiments suggest that this approach has significant potential,
including in the steadily more important context of multi-precision computations
\end{abstract}


\section{Introduction}

We are interested in iterative methods for solving convex quadratic
optimization problems
\begin{equation}\label{quad-prob}
\min_{x \in \mathbb{R}^n} q(x) \eqdef \frac{1}{2}x^T A x - b^Tx
\end{equation}
and large symmetric positive-definite linear systems 
\begin{equation}\label{system}
Ax=b,
\end{equation}
where $A\in\mathbb{R}^{n\times n}$ is symmetric positive definite and $b\in\mathbb{R}^n$.
Such problems are at the centre of efficient methods
in a large variety of domains in applied mathematics, the most prominent
being large-scale numerical nonlinear optimization and
the numerical solution of discretized partial differential equations.
It is thus critical to make the solution
of (\ref{quad-prob}) and (\ref{system}) as efficient as possible.
The cost of most iterative methods for solving these problems
is often dominated by the (potentially many) computations of products
of the form $Ap$ for some vector $p$.
It is therefore of interest to investigate if efficiency gains
may be obtained for this `core' operation.
This is the object of the present paper.

Two different contexts motivate the analysis presented here.
The first is the frequent occurrence of problems involving a hierarchy of model
representations, themselves providing approximations of~$A$ to compute the product $Ap$.
This occurs for instance in discretized applications,
possibly in a multi-resolution framework,
or in inexactly weighted linear and nonlinear least-squares,
where the product itself is obtained by applying an iterative procedure%
\footnote{Our starting point was
  a nonlinear weighted least-squares problem occuring in large-scale data
  assimilation for weather forecasting \cite{GratGuroSimoToin18b}, where the
  inverse of the weighting matrix can not be computed. It thus requires
  the iterative solution of an innermost linear system.}. The second is the
increasing importance of computations in multi-precision arithmetic on the new
generations of high-performance computers
(see \cite{GalaHoro11,PuGalaYangSchaHoro16,cahi18,fahi18,htdh18,HigPraZou19,hipr19} and the many references therein),
in which the use of varying levels of floating point precision
is a key ingredient for obtaining state-of-the-art
energy-efficient computer architectures.
In both cases, using inexact matrix-vector products
(while controlling their inexactness) within the method of conjugate gradients (CG) 
of Hestenes and Stiefel \cite{HestStie52} is a natural option.

Although the use of inexact matrix-vector products in CG
and other Krylov subspace iterative methods has already been investigated
(see for example \cite{SimoSzyl03,vandSlei04,vandSleivanG05,GratToinTshi11,Du.Haber.Karampataki.Szyld.13}),
the proposed analyses typically focus on bounding the Euclidean norm of the residual. This leads to criteria for controlling the inaccuracy of the matrix-vector products  that are somewhat different from the one used here (for instance, see Eq. (4.4) from \cite{Du.Haber.Karampataki.Szyld.13} for the inexact CG,  Eq. (5.8) and (5.9) from  \cite{SimoSzyl03} for the inexact FOM and GMRES). We refer to Section \ref{sec_algo}  for a discussion about the differences between these criteria.  We further refer to  \cite[Section 11]{Simoncini.Szyld.07}, and the many references therein, for a review on inexact Krylov methods.

To the best of our knowledge, none considers the decrease
in the objective function of the associated optimization problem~(\ref{quad-prob}),
which is related to the energy norm of the error in~(\ref{system}). This point of view is, however, important. In optimization, monitoring
the evolution of the nonlinear (and possibly non-convex) objective function or of its model is an obvious concern:
ensuring a fraction of the optimal decrease is, for instance, a standard
convergence argument in trust-region and inexact Newton methods
(see for example\cite[Chapter~6]{ConnGoulToin00}).
In applications arising from elliptic partial differential equations,
several authors have argued that monitoring the energy norm of the error leads
to better termination rules, avoiding under- or over-solving
(see~\cite{ArioDuffRuiz92,Ario04,ArioLoghWath04,Ario13,
              ArioLiesMiedStra13,JiraStraVohr10,PapeLiesStra14}). 
Monitoring the energy norm of the error has been applied for CG with inexact preconditioning, where iterative methods can be embedded for the resolution of the system associated with the computation of the new 
directions.  A stopping criterion in the resolution of the inner systems was suggested in \cite{GolYe99}, that results in a decrease of the energy norm of the residual, and so a decrease of the associated quadratic. However, the accuracy level of the matrix-vector products is fixed.


{\bf Contributions.}
In Section~\ref{rigorous-s} we derive theoretical bounds
on the value of the quadratic in (\ref{quad-prob})
in the presence of inexact matrix-vector products in CG.
We then derive computable estimates required in these theoretical bounds
in Section~\ref{estimates-s}, leading to a practical inexact CG algorithm.
In Section~\ref{numerics-s} we show that very significant efficiency gains
can be obtained in CG by this approach,
both in the case where the accuracy of $Ap$
can be varied continuously and in the case where it is bound
to discrete prescribed levels (as is the case in multi-precision arithmetic).

{\bf Notations.}
Throughout this paper, $\|\cdot\|_2$ denotes the standard Euclidean norm
for vectors and its induced (spectral) norm for matrices.
If~$M$ is symmetric positive definite and $x$ is a vector, $\|x\|_M = \|M^{1/2}x\|_2$.
The dual norm of $\|\cdot\|_M $ with respect
to the Euclidean inner product is~$\|\cdot\|_{M^{-1}}$.
Tr$(M)$ is the trace of the matrix $M$.
If $M$ is symmetric positive definite,
$\lambda_{\min}(M)$ and $\lambda_{\max}(M)$ are its smallest 
and largest eigenvalue, respectively.
$e_i$ is the $i$-th vector of the canonical basis of $\mathbb{R}^n$. 

\section{Analysis of inexact CG}\label{rigorous-s}

We start by stating the CG algorithm with inexact matrix-vector products,
initialized with $x_0=0$.
%

\begin{algorithm}
\caption{Theoretical inexact CG algorithm}\label{alg1}
\begin{algorithmic}
 \State Given symmetric positive definite $A\in\mathbb{R}^{n\times n}$, $b\in\mathbb{R}^n$,
 \State {\it (1)}~~Set $x_0=0$, $\beta_0 = \|b\|_2^2$, $r_0 = -b$, and $p_0=b$ 
  \For{$k=0,1,\ldots,$} 
  \State {\it (2)}~~$c_k =(A+E_k) p_k$
  \State {\it (3)}~~$\alpha_k = \beta_k / p_k^T c_k$
  \State {\it (4)}~~$x_{k+1} = x_k + \alpha_k p_k$
  \State {\it (5)}~~$r_{k+1} = r_k + \alpha_k c_k $ 
  \If{$\|r_{k+1}\|_{A^{-1}}$ is small enough} 
   \State {\it (6)}~~{\bf Stop}
   \EndIf
   \State {\it (7)}~~$\beta_{k+1} = r_{k+1}^T r_{k+1}$
    \State {\it (8)}~~$p_{k+1} = -r_{k+1} + (\beta_{k+1}/\beta_k ) p_k$
  \EndFor
\end{algorithmic}
\end{algorithm}

\noindent

In the above algorithm, the matrix $E_k$ represents a perturbation of $A$
and is the source of inexactness in the matrix-vector product at iteration $k$. 

The name "inexact CG"
 can be viewed as an abuse of language since, due to the error in the matrix-vector product at Step 2 of the algorithm,
  the standard conjugacy of search directions and orthogonality of the residuals is lost. We nevertheless continue to use the designation "inexact CG" because of the  very close similarity between the statements of the method and true CG.

It is known that the residual of the linear system (\ref{system})
provides a handle for monitoring the error in the quadratic $q(x)$ in (\ref{quad-prob}),
provided it is considered in the appropriate norm.
Indeed, if $x_* = A^{-1}b$ is the solution of~(\ref{system})
and $r(x) \eqdef Ax-b$, then
\begin{equation}\label{resA}
\begin{array}{lcl}
\frac{1}{2} \normAi{r(x)}^2
& = & \frac{1}{2} (Ax - b)^TA^{-1}(Ax - b) \\
& = & \frac{1}{2} (x - x_*)^TA(x - x_*) \\
& = & \frac{1}{2} ( x^TAx - 2 x^TAx_* + x_*^TAx_* ) \\
& = & q(x) - q(x_*).
\end{array}
\end{equation}
However, monitoring $\normAi{r(x_k)}$ in inexact CG requires that $r(x_k)$
or a sufficiently good approximation thereof be available,
and that its $A^{-1}$-norm be computed or estimated,
both of which are non-trivial.

When the products $Ap$ are computed inexactly,
the vector $r_k$ recurred in CG is not the same as
the true residual $r(x_k)=Ax_k-b$.
Current literature (see \cite{SimoSzyl03,vandSlei04}) focuses on bounding
the \emph{residual gap} measured in the Euclidean norm, $\|r(x_k)-r_k\|_2$.
Because we are interested in the optimization problem~(\ref{quad-prob}),
we need to bound the residual gap in the $A^{-1}$-norm,
as motivated by (\ref{resA}) and the following lemma.

\begin{lemma}\label{rescombi-l}
Suppose that, at iteration $k$ of the inexact CG algorithm,
\begin{equation}\label{resconds}
\max\Big[ \normAi{r(x_k) - r_k}, \normAi{r_k} \Big]
\leq \frac{\sqrt{\epsilon}}{2} \normAi{b}
\end{equation}
for some $\epsilon>0$.
Then
\begin{equation}\label{quad-decr}
|q(x_k) - q(x_*)| \leq \epsilon |q(x_*)|.
\end{equation}
\end{lemma}

\begin{proof}
  First, evaluating the quadratic $q(x)$	at $x = x_* = A^{-1}b$
	gives a very useful identity, namely that
 \begin{equation}\label{useful}
  |q(x_*)| = -q(x_*) = \frac{1}{2} \normAi{b}^2= \frac{1}{2} \normA{x_*}^2 = \frac{1}{2} |b^Tx_*|. 
  \end{equation}
  Successively using (\ref{resA}), the triangle inequality, (\ref{resconds}), and (\ref{useful}),
	we deduce that
  \begin{equation}\label{qb1}
  \begin{array}{lcl}
    |q(x_k) - q(x_*)|
    &   =  & \frac{1}{2} \normAi{r(x_k)}^2 \\
    & \leq & \frac{1}{2} \left( \normAi{r(x_k) - r_k} + \normAi{r_k}\right)^2 \\
    & \leq & \frac{1}{2} \left( \sqrt{ \epsilon} \normAi{b} \right)^2 \\
    &   =  & \epsilon |q(x_*)|. 
  \end{array} 
  \end{equation} 
\end{proof} 

\noindent
Note that (\ref{quad-decr}) implies
\[
q(x_*) \leq q(x_k) \leq (1-\epsilon) q(x_*).
\]
Thus, if (\ref{resconds}) holds, the quadratic at $x_k$
is within a factor $(1-\epsilon)$ of its minimal value.
Additionally, because $q(x_0)=0$,
\[
|q(x_k)-q(x_0)| \geq (1-\epsilon) |q(x_*)-q(x_0)|.
\]
Thus, if (\ref{resconds}) holds,
the decrease of the quadratic $q(x)$ obtained at $x_k$
is at least $(1-\epsilon)$ times the maximum obtainable decrease.
This is exactly the type of result required to terminate the minimization
in a trust-region context (see \cite[Theorem~6.3.5]{ConnGoulToin00}).

We assume that CG makes   $\Vert r_k\Vert_{A^{-1}}$ small eventually, which it is not guaranteed for the inexact CG. However, if we introduce a reorthogonalization step on the internally-recured residuals, then $r_k$ becomes zero after at most $n$ steps. 
The rest of this paper is devoted to analyzing how to enforce
the part of~(\ref{resconds}) related to the residual gap, that is,
the condition $\normAi{r(x_k) - r_k}\leq \frac{\sqrt{\epsilon}}{2} \normAi{b}$.
Because this last condition measures the residual gap in the $A^{-1}$-norm
(i.e.\ a norm in the dual space), it is natural to use the $A$-norm in the
primal space and the primal-dual matrix norm defined by
\begin{equation}\label{nAiA}
\normAiA{E} \eqdef \sup_{x \ne 0} \displaystyle \frac{\normAi{Ex}}{\normA{x}}
=\| A^{-1/2} E  A^{-1/2}\|_2.
\end{equation}
We will thus use this norm to measure the size of the backward error $E_k$ 
made in the matrix-vector product. Note for future reference that
\begin{equation}\label{nAiAA}
\normAiA{A} = 1
\end{equation}
and 
\begin{equation}\label{nAiAEx}
\normAi{Ex} \leq \normAiA{E}\normA{x}
\end{equation}
for any vector $x$.

%
We first restate, for completeness, a simple result relating
the residual gap to the error matrices $E_j$ 
in inexact CG (see \cite{vandSlei04}).

\begin{lemma}\label{rgapE}
The residual gap in the inexact CG algorithm satisfies
\[
r(x_k) - r_k = - \sum_{j=0}^{k-1}\alpha_j E_j p_j.
\]
\end{lemma}
\begin{proof}
We proceed by induction. 
Observe that $r(x_0)-r_0=0$ and 
\begin{eqnarray*}
  r(x_{1})-r_1
	&=& (Ax_{1} - b) - r_{1} \\
  &=& (\alpha_0 Ap_0 - b ) - (r_0 + \alpha_0 c_0) \ = \ - \alpha_0 E_0 p_0.
\end{eqnarray*}
Suppose now that the result is true for iterations $j=0, \ldots, k$.
From the recurrences for $x_{k+1}$ and $r_{k+1}$ we then have that 
\begin{eqnarray*}
  r(x_{k+1})-r_{k+1}
	&=& (Ax_{k+1} - b) - r_{k+1} \\
  &=& (Ax_k+\alpha_k Ap_k - b ) -r_k + r_k -r_{k+1} \\
  &=& r(x_k) +\alpha_k Ap_k - r_k - \alpha_k(A+E_k)p_k\\
  &=& r(x_k) - r_k - \alpha_k E_k p_k,
\end{eqnarray*}
from which the result follows.
\end{proof} 

\noindent
We are now in a position the derive suitable bounds on the error matrices.

\begin{theorem}\label{CG-bounds}
Let $\epsilon>0$ and let $\phi \in \mathbb{R}^k$ be a positive vector satisfying
\begin{equation}\label{sumphiinv}
\sum_{j=1}^{k}\frac{1}{\phi_j} \leq 1.
\end{equation}
Suppose furthermore that 
\begin{equation}\label{nEj-CG}
\normAiA{E_j}
\leq \omega_j
\eqdef \frac{\sqrt{\epsilon} \,\normAi{b} \normA{p_j}}
         {2\,\phi_{j+1} \|r_j\|_2^2 + \sqrt{\epsilon} \,\normAi{b} \normA{p_j} }
\end{equation}												
for all $j\in\{ 0, \ldots, k-1 \}$.
Then
\[
\normAi{r(x_k)-r_k} \leq \frac{\sqrt{\epsilon}}{2}\normAi{b}.
\]
If additionally
\begin{equation}\label{condr}
\normAi{r_k}\leq \frac{\sqrt{\epsilon}}{2}\normAi{b},
\end{equation}
then (\ref{resconds}) and (\ref{quad-decr}) both hold. 
\end{theorem}

\begin{proof} First note that (\ref{nEj-CG}) ensures that $\omega_j \in (0,1)$.
  Lemma~\ref{rgapE}, the triangle inequality, and~(\ref{nAiAEx}) imply that
  \begin{equation}\label{cg1}
	\begin{array}{lcl}
  \normAi{r(x_k) - r_k}
   &\leq &\displaystyle \sum_{j=0}^{k-1} \normAi{\alpha_j E_j p_j} \\
   &\leq &\displaystyle \sum_{j=0}^{k-1} |\alpha_j| \normAiA{E_j} \normA{p_j}. 
	\end{array}
  \end{equation}
  Now, using~(\ref{nAiAEx}) again and~(\ref{nEj-CG}),
  \[
  \alpha_j
    = \frac{\|r_j\|_2^2}{p_j^T(A+E_j)p_j} 
    \leq \frac{\|r_j\|_2^2}{p_j^TAp_j - \normAiA{E_j}\normA{p_j}^2} 
	  \leq \frac{\|r_j\|_2^2}{(1-\omega_j\,)\normA{p_j}^2}.  
  \]
  Substituting this bound in (\ref{cg1}) and using (\ref{nEj-CG}) again, we obtain 
\begin{equation}\label{cg2}
  \normAi{r(x_k) - r_k}
  \leq \sum_{j=0}^{k-1} \frac{\omega_j}{1-\omega_j} \frac{\|r_j\|_2^2}{\normA{p_j}}.
  \end{equation}
  But the definition of $\omega_j$ in (\ref{nEj-CG}) gives 
  \[
  \displaystyle \frac{\omega_j}{1-\omega_j}
   =\displaystyle
	     \frac{\sqrt{\epsilon}\,\normAi{b}\normA{p_j}}{2\,\phi_{j+1}\|r_j\|_2^2},
  \]
  so that (\ref{cg2}) becomes
  \begin{equation}\label{CG-bound}
	\begin{array}{lcl}
  \normAi{r(x_k) - r_k}
   &\leq& \displaystyle
	     \sum_{j=0}^{k-1} \frac{\sqrt{\epsilon}\,\normAi{b}\normA{p_j}}{2\,\phi_{j+1}\|r_j\|_2^2} 
	                      \frac{\|r_j\|_2^2}{\normA{p_j}} \\
   &=& \displaystyle
	     \frac{\sqrt{\epsilon}}{2}\normAi{b} \sum_{j=0}^{k-1} \frac{1}{\phi_{j+1}} 
	      \ \leq \ \frac{\sqrt{\epsilon}}{2}\normAi{b}.
	\end{array}
  \end{equation}
	The result then follows from Lemma~\ref{rescombi-l}.
\end{proof} 

\noindent
Observe that (\ref{nEj-CG}) allows a perturbation in $A$
whose norm $\normAiA{E_j}$ depends on the ratio $\displaystyle{\frac{\Vert r_j\Vert^2_2}{\Vert p_j\Vert_A}}$.
In particular, when this ratio decreases as the iterations proceed
(although the decrease of $\Vert r_j\Vert_2$ is not generally monotonic in CG)
the error in the matrix-vector products is permitted to grow.

Some additional comments on Theorem~\ref{CG-bounds} are in order at this stage.
\begin{enumerate}
\item (\ref{nEj-CG}) assumes that the primal-dual norm is the natural norm
  for measuring the size of the error matrices $E_j$.
	While this may be true in certain applications,  
  $\normAiA{E_j}$ may be difficult to compute or estimate in practice.
\item Even discounting that potential difficulty, 
  verifying conditions~(\ref{nEj-CG}) and~(\ref{condr}) remains impractical,
	as the quantities $\normA{p_j}$, $\normAi{b}$, and $\normAi{r_k}$
	are not readily available in the course of the inexact CG algorithm.
\item The $\phi_j$ appearing in (\ref{nEj-CG}) may be used as
  part of a global ``error management strategy''.
	A simple choice that obviously satisfies (\ref{sumphiinv})
	is to define $\phi_j = k_{\max}$ for all $j$,
	where $k_{\max}$ is the maximum allowable number of iterations.
	In fact, the $\phi_j$ can be used to further advantage.
\end{enumerate}

\noindent
We shall address the above issues in the following section. 

\section{A practical inexact CG algorithm}\label{estimates-s}

\subsection{Managing the inaccuracy budget}\label{inacc_budget-s}

An important ingredient of a practical inexact CG algorithm
is the choice of the $\phi_j$ in (\ref{nEj-CG}).
Note that (\ref{sumphiinv}) constrains the $\phi_j$
over \emph{all iterations until termination}.
As mentioned earlier, choosing $\phi_j=k_{\max}$
is adequate but often suboptimal.
Indeed, it is possible to adjust the $\phi_j$ adaptively,
in particular when the inaccuracy of the product $Ap$ cannot be varied continuously
but is bound to prescribed levels
(for example, different levels of floating point precision).

Suppose that for a given $\phi_{j+1}$, 
an inexactness $\omega_j(\phi_{j+1})$ in the matrix-vector product at step $j$
is allowed by Theorem~\ref{CG-bounds}, but the actual error~$E_j$
satisfies $\normAiA{E_j} = \hat{\omega}_j < \omega_j$.
This implies that a larger value $\hat\phi_{j+1}$
could have been used instead of $\phi_{j+1}$.
Solving for $\hat{\phi}_{j+1}$ in the linear equation
$\hat{\omega}_j = \omega_j(\hat{\phi}_{j+1})$ in~(\ref{nEj-CG}) gives 
\begin{equation}\label{hatphi}
\hat{\phi}_{j+1}
     = \frac{1-\hat{\omega}_j}{\hat{\omega}_j}
		   \frac{\sqrt\epsilon\,\normAi{b} \normA{p_j}}{2\|r_j\|_2^2} > \phi_{j+1}.
\end{equation}
We may then distribute the unused inaccuracy $1 - \sum_{p=1}^{j+1}\hat{\phi}_p^{-1}$
evenly in the remaining $k_{\max}-j-1$ iterations,
by setting
\[
\phi_i = \frac{k_{\max}-j-1}{1 - \sum_{p=1}^{j+1}\hat{\phi}_p^{-1}},
\;\;\;\;
i = j+2, \ldots, k_{\max}.
\]
This leads to smaller values of $\phi_i$ 
(and therefore larger allowable errors)
in subsequent iterations.
The updated $\phi_i$ still satisfies~(\ref{sumphiinv}),
as shown below:
\begin{align*}
\sum_{i=1}^{k_{\max}} \phi_i^{-1}
 &= \sum_{i=1}^{j+1} \phi_i^{-1}
    + \sum_{i=j+2}^{k_{\max}} \phi_i^{-1} \\
 &= \sum_{i=1}^{j+1} \phi_i^{-1}
      + (k_{\max}-j-1) \frac{ 1 - \sum_{p=1}^{j+1}\hat{\phi}_p^{-1} }{ k_{\max}-j-1 } \\
 &= \sum_{i=1}^{j+1} \phi_i^{-1} + 1 - \sum_{i=1}^{j+1}\hat{\phi}_i^{-1} \ < \ 1,
\end{align*}
since $\hat{\phi}_i>\phi_i$.
In practice, this allows maintaining only single
running values for $\phi_{j+1}$ and
\[
\Phi_j \eqdef 1-\sum_{p=1}^j\hat{\phi}_p^{-1}
\]
for $j$ ranging from $0$ to $k_{\max}-1$.

\subsection{Computable estimates of $\normAiA{E_j}$ and $\normA{p_j}$}

We now attempt to estimate the quantities required by Theorem~\ref{CG-bounds}
that are unavailable in the inexact CG algorithm.

We first consider that $\normAiA{E_j}$ in (\ref{nEj-CG})
may not be available from the application context
and note that, from (\ref{nAiA}),
\begin{equation}\label{app-nEj}
\normAiA{E_j}
= \| A^{-1/2} E_j A^{-1/2} \|_2
\leq \lambda_{\min}(A)^{-1} \|E_j\|_2,
\end{equation}
so that a bound on $\|E_j\|_2$ can be used provided one knows
(an approximation of) the smallest eigenvalue of $A$.
To estimate $\normA{p_j}$, we can use the fact that
\[
\lambda_{\min}(A)^{1/2} \|p_j\|_2 \leq \normA{p_j} \leq\lambda_{\max}(A)^{1/2}  \|p_j\|_2.
\]
However, for ill-conditioned problem, the above bounds are likely to be very loose.
Another approach is to choose
\begin{equation}\label{app-pj}
\normA{p_j} \approx \sqrt{\frac{1}{n}{\rm Tr}(A)} \|p_j\|_2.
\end{equation}
This can be justified by the fact that
each side of the above expression would have the same mean squared value 
if the entries of $p_j$ were independent standard normal variables.

\subsection{A computable estimate of $\normAi{b}$}

Finding an estimate of $\normAi{b}$ is more difficult,
since this quantity is related to the value of the quadratic $q(x)$
\emph{at the solution}~$x_*$, see~(\ref{useful}).
Note from~(\ref{resA}) that
\[
|q(x_*)| \leq |q(x_k)| + \frac{1}{2} \normAi{r(x_k)}^2.
\]
To our knowledge, the best available approximation is
the absolute value of the quadratic at the current iterate,
and thus we choose
\begin{equation}\label{nAib-est}
\normAi{b} = \sqrt{2 |q(x_*)|}
           \approx \sqrt{2 |q(x_k)|}. 
\end{equation}
If there is no error in the products $Ap_j$ in CG,
assuming exact arithmetic, $r(x_k)$ is orthogonal to~$x_k$
and it follows that 
\[
q(x_k) = \frac{1}{2} x_k^TAx_k - b^Tx_k =
         \frac{1}{2} x_k^T(Ax_k-b)  - \frac{1}{2} b^T x_k = -\frac{1}{2} b^Tx_k.
\]
In the presence of errors in the matrix-vector products,
the above property may no longer hold,
and it is of interest to analyze how much
\begin{equation}\label{qk}
q_k \eqdef -\frac{1}{2} b^Tx_k
\end{equation}
differs from $q(x_k)$.
This is also important if the decrease in the quadratic objective
function is used for other purposes, as is the case, for instance,
in trust-region methods, where it is a key ingredient in the management
of the trust-region radius (see \cite[Chapter~6]{ConnGoulToin00}).
To this end, we first prove the following backward error property.

\begin{lemma}\label{backw-slope}
Let $A\in\mathbb{R}^{n\times n}$ be symmetric positive definite
and $b\in\mathbb{R}^n$. For any nonzero $x\in\mathbb{R}^n$, with $r(x)=Ax-b$,
\begin{equation}\label{slope-pert}
\min_{E} \left\{ \frac{\normAiA{E}}{\normAiA{A}} \mid x^T(A+E)x= x^Tb \right\}
= \frac{|x^Tr(x)|}{\normA{x}^2}.
\end{equation}
\end{lemma}

\begin{proof}
  Using (\ref{nAiA}) and (\ref{nAiAA}), we rewrite the optimization problem
  as
  \[
  \min_E \left\{ \|A^{-1/2}EA^{-1/2}\|_2 \mid x^TEx = x^T(b - Ax) \right\}
  \]
  which, since $A$ is positive definite and setting $y = A^{1/2}x$,
	is itself equivalent to 
  \[
  \min_E  \left\{ \|A^{-1/2}EA^{-1/2}\|_2\mid y^TA^{-1/2}EA^{-1/2}y = -x^Tr(x) \right\}.
  \]
  But for any nonzero vector $y$ and scalar $\gamma$, 
  \[
  \min_M  \left\{ \|M\|_2 \mid y^TMy = \gamma \right\}
  = \frac{\gamma}{\|y\|_2^2}
  \]
  and the minimum is attained by $M = \gamma yy^T/\|y\|_2^4$ with $\|M\|_2 = \gamma/\|y\|_2^2$.
	Thus the minimum in (\ref{slope-pert}) is
  \[
  \normAiA{E} =\|A^{-1/2}EA^{-1/2}\|_2 = \frac{|x^Tr(x)|}{\normA{x}^2}.
  \]  
\end{proof}

\noindent
We can now show that $q_k$ remains close to $q(x_k)$ despite 
the inexact matrix-vector products.

\begin{theorem}\label{errqval-th}
  Let $x_k$ be the result of applying the inexact CG algorithm
	and let $q_k = -\frac{1}{2} b^Tx_k$.
	If the inexactness of the matrix-vector products is controlled as in (\ref{nEj-CG}),
	and (\ref{condr}) holds, then
  \begin{equation}\label{errqval}
  \frac{|q(x_k) - q_k|}{|q(x_*)|} \leq \frac{\sqrt{\epsilon}(1+\sqrt{\epsilon})}{2}.
  \end{equation}
\end{theorem}

\begin{proof}
  We deduce from Lemma~\ref{backw-slope} that there exists a quadratic 
   \begin{equation}\label{hatq-def}
  \hatq(x_k) = \frac{1}{2} x_k^T(A+E)x_k - b^Tx_k  
  \end{equation}
  such that by construction $\hatq(x_k) = -\frac{1}{2} b^Tx_k = q_k$ and
   \begin{equation}\label{weak}
  \begin{array}{lcl}
  |q(x_k) - q_k|
    & =  & |q(x_k) - \hatq(x_k)| \\*[2ex]
    & =  & \frac{1}{2} |x_k^TEx_k| \\*[2ex]
    &\leq& \frac{1}{2} \normAiA{E}\normA{x_k}^2 \\*[2ex]
    & =  & \frac{1}{2} |x_k^Tr(x_k)|\\*[2ex]
    &\leq& \frac{1}{2} \normAi{r(x_k)}\normA{x_k},
  \end{array}
  \end{equation}
  where we used (\ref{hatq-def}), Lemma~\ref{backw-slope},
	and the Cauchy-Schwarz inequality.
  But by Theorem~\ref{CG-bounds} and (\ref{useful}),
  \[
  \begin{array}{lcl}
  \normA{x_k}
    &\leq& \normA{x_*} + \normA{x_k-x_*} \\
    & =  & \normA{x_*} + \normAi{r(x_k)} \\
    &\leq& \normA{x_*} + \normAi{r(x_k)-r_k} + \normAi{r_k}\\
    &\leq& ( 1 + \sqrt{\epsilon} ) \normA{x_*}
  \end{array}
  \]
  and hence, using (\ref{condr}) and (\ref{useful}) again,
  \[
   |q(x_k) - q_k|
   \leq \frac{1}{2} \normAi{r(x_k)} (1+\sqrt{\epsilon}) \normA{x_*}
   \leq \frac{1}{2}\sqrt{\epsilon}(1+\sqrt{\epsilon})|q(x_*)|.
  \]
  \end{proof}
	
\noindent
The above theorem shows that $q_k$ in (\ref{qk}) remains close to $q(x_k)$.
The bound is considerably weaker than (\ref{quad-decr}),
but it is likely to be pessimistic
as we have not taken into account in (\ref{weak}) the fact that the
angle between $x_k$ and $r(x_k)$ is expected to be small.
(For CG in exact arithmetic with exact matrix-vector products, $x_k^Tr_k=0$.)
This will be numerically confirmed in Section~\ref{numerics-s}.

From~(\ref{nAib-est}) and Theorem \ref{errqval-th},
we use the approximation
 \begin{equation}\label{app-nAib}
	\normAi{b} \approx \sqrt{2|q_k|},
\end{equation}
for $k=1,\ldots,k_{\max}$ 
in our practical inexact CG algorithm.
For $k=0$, because $x_0=0$, we use the rougher approximation
$\normAi{b} \approx \|b\|_2 / \sqrt{\lambda_{\max(A)}}$. This conservative estimate is likely to result with the computation of the first matrix-vector product in double precision. Different values could be naturally chosen in $\displaystyle{\left[ \frac{\Vert b\Vert _2}{\sqrt{\lambda_{max(A)}}} , \frac{\Vert b\Vert _2}{\sqrt{\lambda_{min(A)}}}\right]}$.

\subsection{A computable estimate of $\normAi{r_k}$}

It is also necessary to estimate $\normAi{r_k}$,
in order to perform the termination test~(\ref{condr}).
Estimating the energy norm of the error in CG is a well-studied problem
(see for example~\cite{StraTich02,MeurStra06,JiraTitl10,ArioGrat12,LiesStra13}).
We follow these ideas, ignoring pathological convergence instances,
and estimate
\[
\frac{1}{2} \normAi{r_k}^2 \approx q(x_k) - q(x_*) \approx q(x_{k-d})-q(x_k)
\]
where $d$ is a small integer ($d=10$ in our case).
Using this estimate in (\ref{condr}) leads to 
\[
q(x_{k-d})-q(x_k) \leq \frac{1}{4} \epsilon |q(x_k)|,
\]
which, using (\ref{app-nAib}), is itself approximated by
 \begin{equation}\label{delay-termination}
q_{k-d}-q_k \leq \frac{1}{4} \epsilon |q_k|.
\end{equation}

\subsection{Resulting practical algorithm}
\label{sec_algo}

We now consolidate our approximations in order to obtain
a practical version of the inexact CG algorithm 
that relies only on computable quantities.
We note that our definitions of these quantities nevertheless
requires the user to provide a (potentially very rough, see Section~\ref{numerics-s})
approximation of the smallest and largest eigenvalues of $A$.

Making use of (\ref{app-nEj}),~(\ref{app-pj}), and~(\ref{app-nAib}),
we suggest to approximate (\ref{nEj-CG}) by
 \begin{equation}\label{app-nEj-CG}
\frac{\|E_j\|_2}{\lambda_{\min}(A)}
 \leq \frac{\sqrt{\epsilon}\,\sqrt{|q_j|}\,\sqrt{{\rm Tr}(A)}\,\|p_j\|_2}
     {\sqrt{2n}\,\phi_{j+1}\|r_j\|_2^2 
		   + \sqrt{\epsilon}\,\sqrt{|q_j|}\,\sqrt{{\rm Tr}(A)}\,\|p_j\|_2},
\end{equation}
for $j=1,\ldots\,k-1$.
The formula is similar for $j=0$, with $\sqrt{|q_0|}$
replaced by $\sqrt{2}\|b\|_2/\sqrt{\lambda_{\max}(A)}$.
We also replace the termination test (\ref{condr})
by its practical version (\ref{delay-termination}).

Our practical bound  differs from the one suggested by \cite{Du.Haber.Karampataki.Szyld.13} and defined as  $$ \Vert E_j\Vert_2\leq \frac{\sigma_{\min}(A)}{2} \min(1,\frac{\sqrt{\epsilon} \Vert p_j\Vert_2}{m\Vert r_j\Vert_2^2}), \quad \forall j=0,1,\cdots,m-1. $$
However, we note that the two bounds depend in the same manner of the smallest eigenvalue of $A$ and  the ratio $\displaystyle{\frac{\Vert r_j\Vert_2^2}{\Vert p_j\Vert_2}}$ when it is large. The decrease in $\displaystyle{\frac{\Vert r_j\Vert_2^2}{\Vert p_j\Vert_2}}$ allows the increase of the bounds of the error in the matrix-vector products, while the smallest eigenvalue of $A$ may constrain the bound to remain tiny.

\begin{algorithm}
\caption{Practical inexact CG algorithm}\label{practical-CG}
\begin{algorithmic}
 \State Given symmetric positive definite $A\in\mathbb{R}^{n\times n}$, $b\in\mathbb{R}^n$, $\epsilon$, $k_{\max}$, $\lambda_{\min}$, $\lambda_{\max}$, and reorth.
 \State {\it (1)}~~Set $x_0=0$, $r_0 = -b$, $q_0= 0$, $\beta_0= \|b\|_2^2$, $u_1= b / \beta_0$, $p_0=b$, $\phi_0 = k_{\max}$, and $\Phi_0 = 1$
  \For{$k=0,1,\ldots,k_{\max}$} 
  \State {\it (2)}~~Determine $\omega_k$ defined by the right hand side of Eq. (\ref{app-nEj-CG})
  \State {\it (3)}~~Compute the product $c_k =(A+E_k)p_k$ with $\normAiA{E_k} \leq \omega_k$, also returning $\hat{\omega}_k$
  \State {\it (4)}~~$\alpha_k = \beta_k / p_k^T c_k$
  \State {\it (5)}~~$x_{k+1} = x_k + \alpha_k p_k$ 
  \State {\it (6)}~~ $q_{k+1} = \frac{1}{2} b^Tx_{k+1}$
  \If{($q_{k+1-d}-q_{k+1}) \leq \frac{1}{4} \epsilon |q_{k+1}|$ } 
   \State {\it (7)}~~{\bf Stop}
   \EndIf
   \State {\it (8)}~~Compute $\hat{\phi}_k$ from $\hat{\omega}_k$ using the first part of Eq. (\ref{hatphi})
    \State {\it (9)}~~$\Phi_{k+1} = \Phi_k - \hat{\phi}_k^{-1}$
   \If{$k<k_{\max}$}
    \State {\it (10)}~~$\phi_{k+1} =(k_{\max}-k)/\Phi_{k+1}$
     \Else
    \State {\it (11)}~~ $\phi_{k+1} =\phi_{k}$
     \EndIf
    \State {\it (12)}~~$r_{k+1} = r_k + \alpha_k c_k $
    \If{(reorth)}
    \For{$i=1,\ldots,k$}    
     \State {\it (13)}~~$r_{k+1} = r_{k+1} - (u_i^Tr_{k+1}) u_i$
    \EndFor
    \State  {\it (14)}~~$\beta_{k+1} = r_{k+1}^T r_{k+1}$
     \State  {\it (15)}~~$u_{k+1} = r_{k+1} / \sqrt{\beta_{k+1}}$
     \Else
     \State  {\it (16)}~~$\beta_{k+1} = r_{k+1}^T r_{k+1} $
    \EndIf
     \State  {\it (17)}~~$p_{k+1} = -r_{k+1} + (\beta_{k+1}/\beta_k)p_k$
  \EndFor
\end{algorithmic}
\end{algorithm}

We include the option of reorthogonalization 
in our practical inexact CG method. It is well known that the residual vectors, though theoretically orthogonal,
quickly lose their orthogonality in finite precision computations.
The reorthogonalization step is applied to the sequence of the internally-recurred residuals using the MGS algorithm prior to the computation of the new search direction. It means that the $A$-norm of the error is minimized over some expanding subspaces. To the limits of the machine precision, the iterative process is expected to converge in 
 $n$ steps at most, and guarantees that  $\Vert r_k\Vert_{A^{-1}}$ will become  lower than the prescribed tolerance.

However, it is well beyond the scope of this paper to analyze the effects
of finite precision and the loss of orthogonality on the convergence
of CG with inexact matrix-vector products,
but we do report on some numerical experiments
with and without reorthogonalization in the following section.

\section{Numerical experiments}\label{numerics-s}

We first provide some figures to illustrate 
typical behaviour of the theoretical and practical inexact CG algorithms.

\subsection{Continuously varying precision}

In our first example, $A$ is a $100\times100$ diagonal matrix
with entries logarithmically equally-spaced
between $1$ and $10^{-4}$, and $b=A[1,\dots,1]^T$.

\begin{figure}[htb!]

\begin{center}
\begin{tabular}{cc}
{\footnotesize (a) double precision CG}     &
{\footnotesize (b) theoretical inexact CG}  \\
\includegraphics[width=3in]{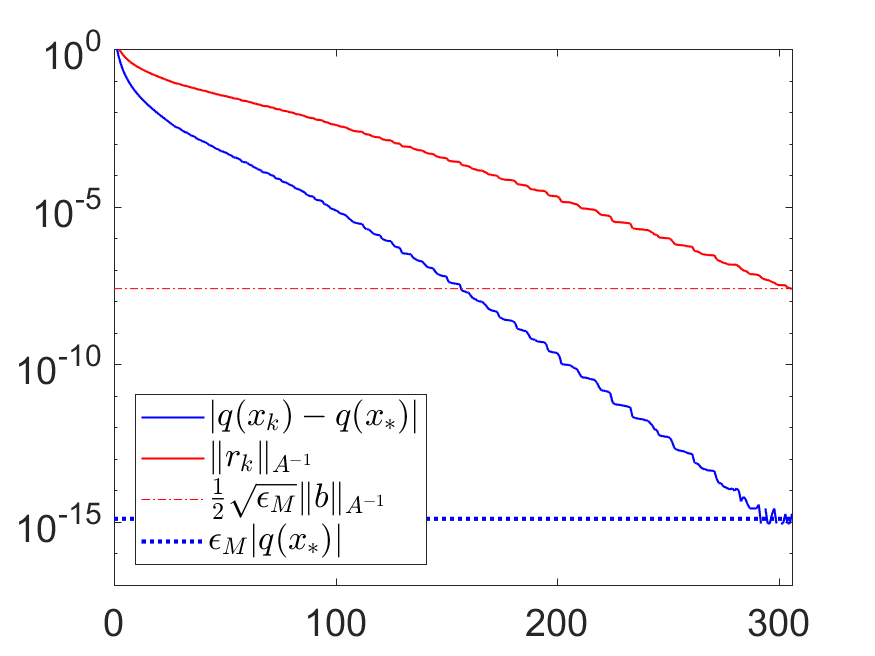} &
\includegraphics[width=3in]{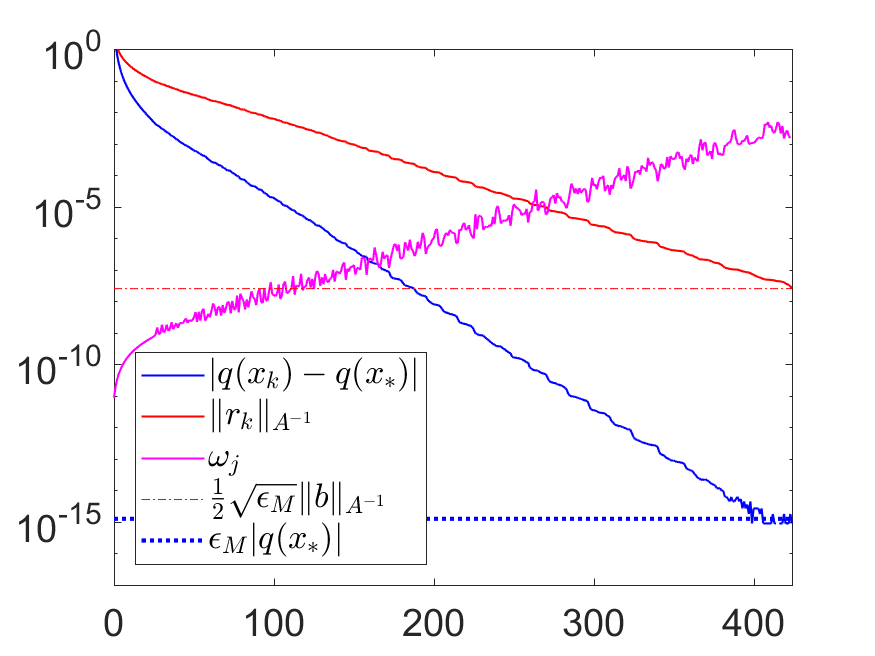} \\
{\footnotesize (c) practical inexact CG}     &
{\footnotesize (d) with reorthogonalization} \\
\includegraphics[width=3in]{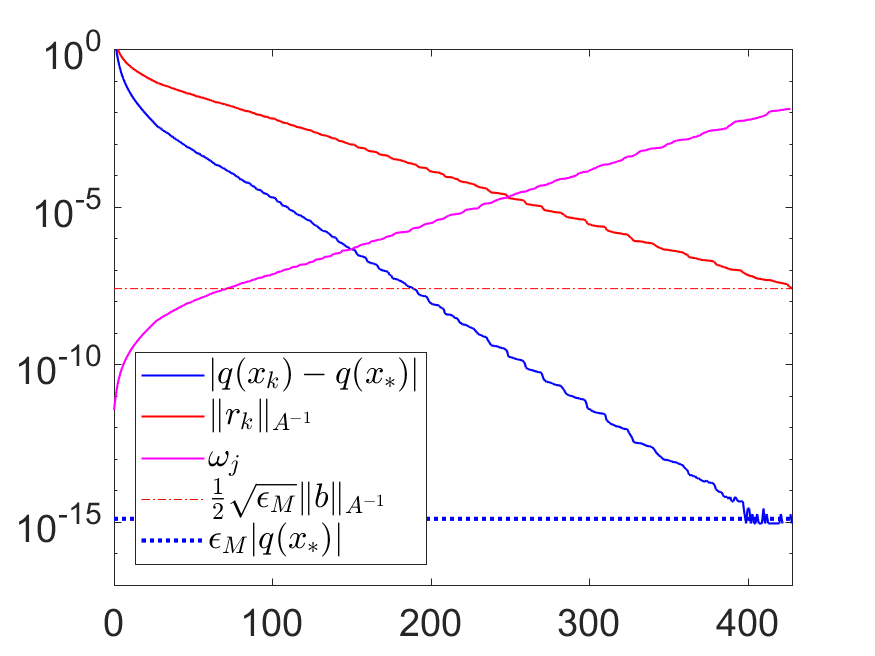} &
\includegraphics[width=3in]{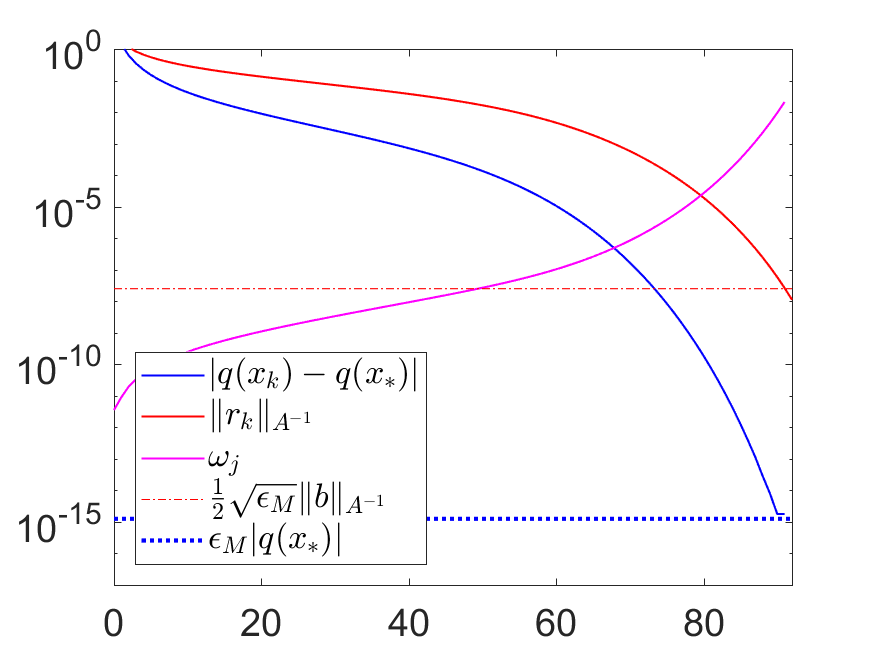} 
\end{tabular}
\caption{CG applied to $\mathrm{diag(logspace(-4,0,100))}$}
\label{eg1}
\end{center}

\end{figure}

\noindent
We plot both the decrease in the quadratic 
and the $A^{-1}$ norm of the residual.
We terminate the iteration when~(\ref{condr}) holds
with $\epsilon=\epsilon_M=2^{-52}$, the IEEE double machine precision. 
Lemma~\ref{rescombi-l} then ensures that~(\ref{quad-decr}) holds.

In Figure~\ref{eg1}(a) the matrix-vector products are performed
in full double precision arithmetic.
In Figure~\ref{eg1}(b) they are performed inexactly,
with a random perturbation $E$ satisfying
the theoretical condition~(\ref{nEj-CG}).
In Figure~\ref{eg1}(c) the matrix-vector products
are once again performed inexactly,
but with a random perturbation $E$ satisfying
the practical condition~(\ref{app-nEj-CG}).
Figure~\ref{eg1}(d) is the same as Figure~\ref{eg1}(c),
except CG is performed with full reorthogonalization. 

\begin{figure}[htb!]

\begin{center}
\begin{tabular}{cc}
{\footnotesize (a) double precision CG}     &
{\footnotesize (b) theoretical inexact CG}  \\
\includegraphics[width=3in]{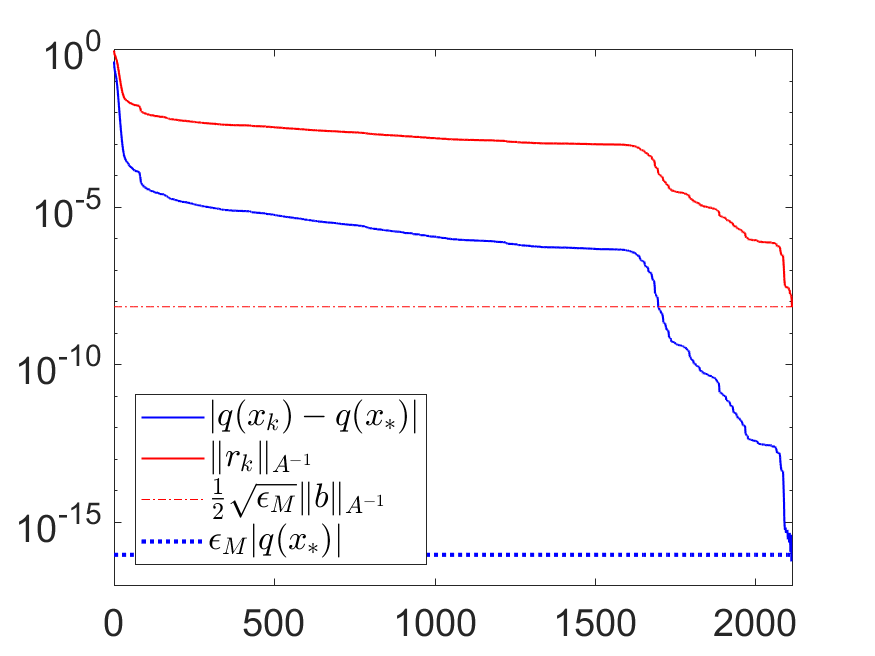} &
\includegraphics[width=3in]{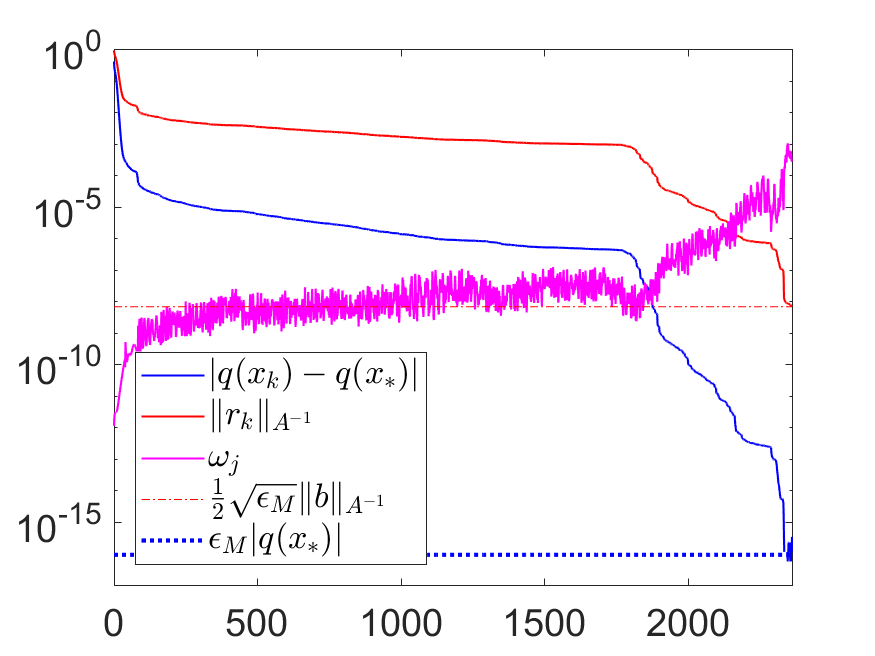} \\ 
{\footnotesize (c) practical inexact CG}     &
{\footnotesize (d) with reorthogonalization} \\
\includegraphics[width=3in]{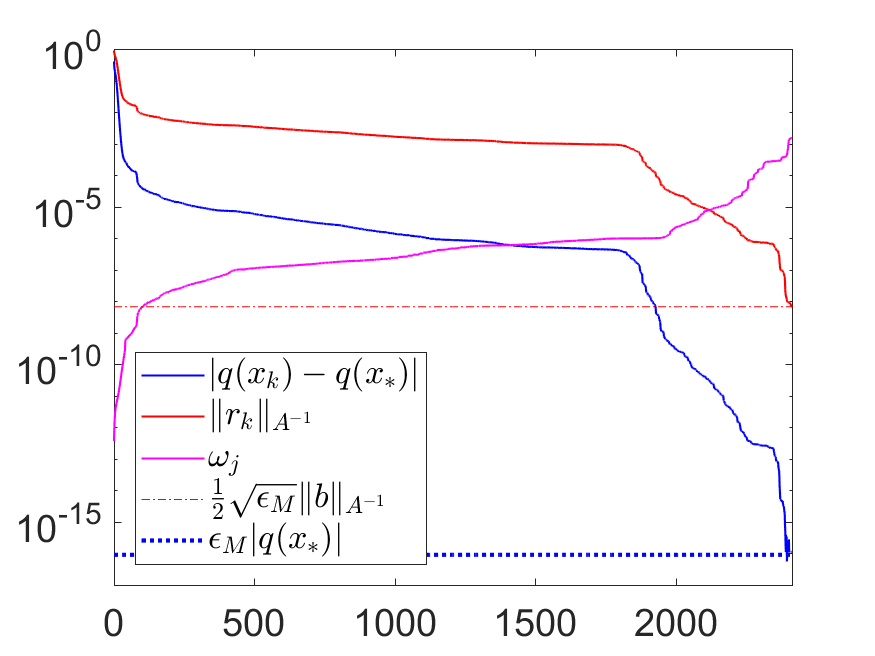} &
\includegraphics[width=3in]{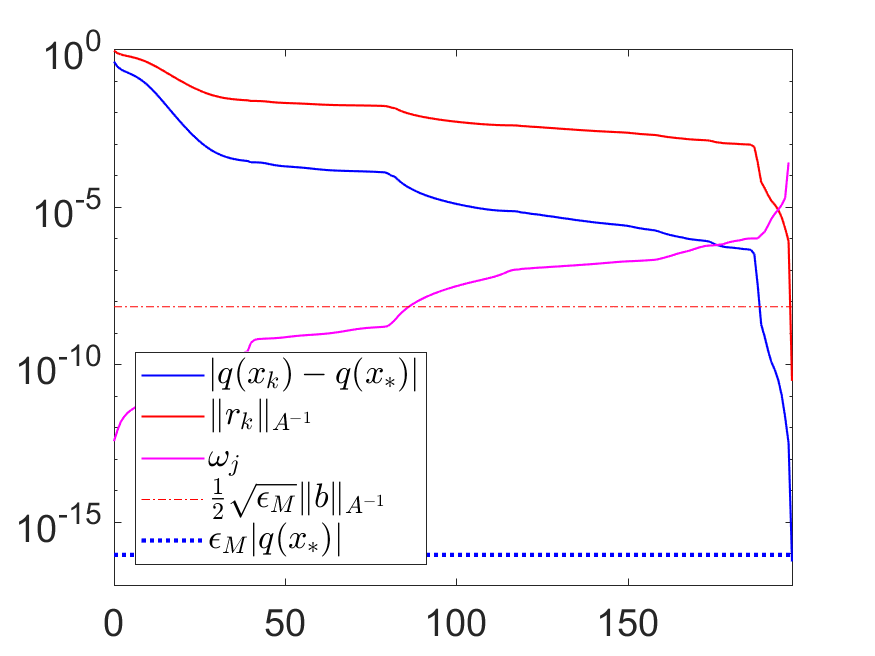} 
\end{tabular}
\caption{CG applied to $\mathrm{nos1.mat}$}
\label{eg2}
\end{center}

\end{figure}

In our second example, we repeat the same tests
with $A$ the matrix $\mathrm{nos1.mat}$ from the Matrix Market,
scaled to have Euclidean norm 1.
This is a $237\times237$ matrix
with condition number approximately $10^8$.
Results are shown in Figure~\ref{eg2}.

In both examples, performing inexact matrix-vector products 
may lead to delays in convergence compared
to double precision CG.
The quantity $\omega_j$ in (\ref{nEj-CG})
can be quite oscillatory (see figures~(b)).
When it is estimated as in (\ref{app-nEj-CG}), however,
it seems to increase monotonically (see figures~(c)).
In these examples, despite the heuristic nature of the estimates,
the practical criterion (\ref{app-nEj-CG}) works just as well as
the theoretical criterion (\ref{nEj-CG}) (compare figures~(b) and~(c)).
Without reorthogonalization,
CG may require (much) more than $n$ iterations to converge
to the required tolerance,
especially in the second ill-conditioned problem.
In both these examples, however, 
reorthogonalization ensures convergence of
the practical inexact CG algorithm
in fewer than $n$ iterations (see figures~(d)).

\subsection{Discrete precision levels}

In the next example, we suppose that only IEEE double, single and half precisions are available. We will refer to IEEE precisions when using the expressions double, simple or half precisions until the end of the manuscript. We switch to matrix-vector products in a lower precision
once the error in the lower precision satisfies
the practical condition~(\ref{app-nEj-CG}).
We modify the $\phi_j$ adaptively as explained in Section~\ref{inacc_budget-s}. In practice, the use of the single and half precisions is likely to trigger under- and overflows both when converting numbers in the targeted format and computing the matrix-vector products.  However, \cite{HigPraZou19} suggested an algorithm that prevent overflows when converting a matrix from double or single precision to half precision. Nevertheless if occurence of such issues, the precision and the inaccuracy budget will have to be adapted accordingly (computation in simple or double precision). This would certainly prevent the use of the IEEE half precision depending on the application. We did not implement such strategy nor assess the occurence of this issue in the following  numerical experiments due to the use of emulated accuracy.

\begin{figure}[ht]

\begin{center}
\begin{tabular}{cc}
{\footnotesize (a) $\mathrm{diag(logspace(-4,0,100))}$}   &
{\footnotesize (b) $\mathrm{nos1.mat}$}                   \\
\includegraphics[width=3in]{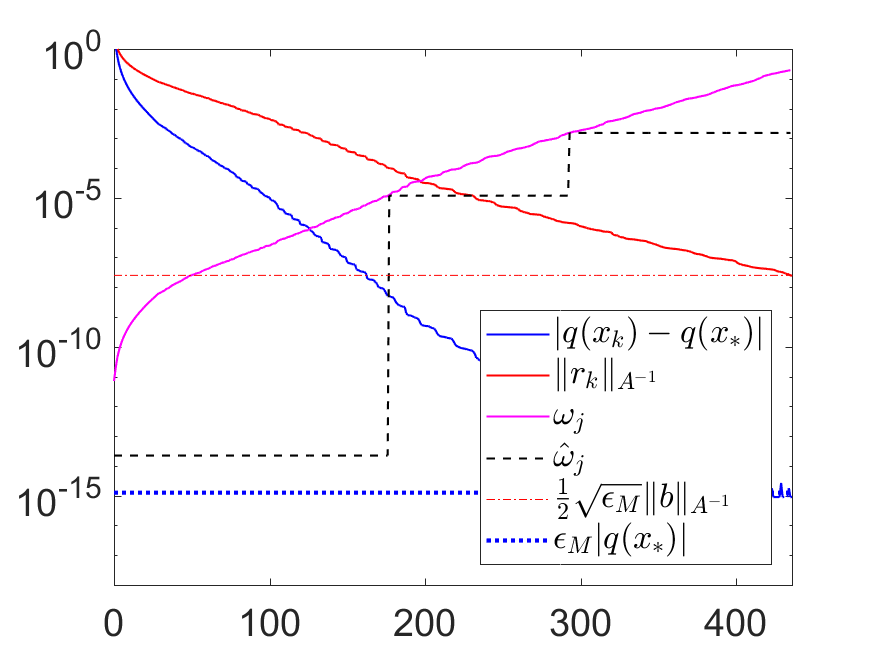} &
\includegraphics[width=3in]{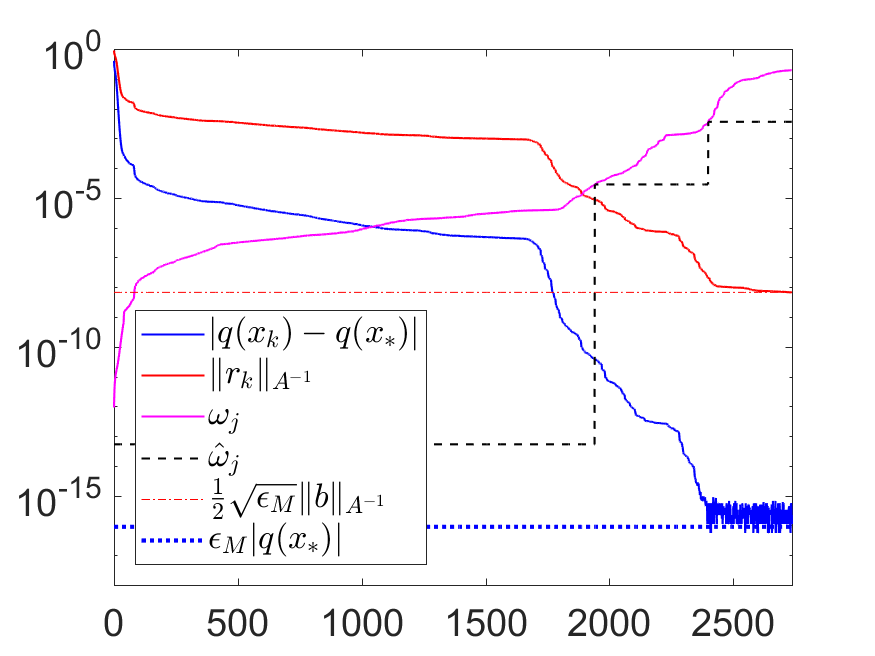} 
\end{tabular}
\caption{inexact CG in discrete precision levels}
\label{eg3}
\end{center}

\end{figure}

Results are shown in Figure~\ref{eg3}.
The convergence behaviour is similar
to that of the continuously-varying precision case.

\subsection{Performance comparisons}

We consider 4 algorithms:

\vspace*{1mm}
\begin{tabular}{lp{5in}}
CG:   & the standard CG algorithm with products computed in double precision, \\
CGR:  & the standard CG algorithm with reorthogonalization and products computed in double precision \\
iCG:  & the inexact CG algorithm without reorthogonalization,\\
iCGR: & the inexact CG algorithm with reorthogonalization (computed in double precision). 
\end{tabular}
\vspace*{0.25cm}

The reorthogonalization step introduced in CGR and iCGR corresponds to the
Gram-Schmidt process with the standard inner product. In both algorithms, it
is applied to the sequence of the recurred residual $\{r_k\}$ and is computed
in double precision. iCGR differs from the flexible CG algorithm \cite{Not00}
in the sense that the Gram-Schmidt process is not used to prevent the loss of
the $A$-conjugacy of the search directions $\{p_k\}$, but rather a loss of
orthogonality of the recurred residual $\{r_k\}$. Therefore, the practical
variant of iCGR does not require computing the matrix-vector products with
$A$ in double precision, at variance with the flexible CG algorithm.

In order compare the performance of these algorithms,
we assume that the computational cost, and so the energy cost, is dominated
by the matrix-vector products, and use the following two models of energy cost. 
For the continuously-varying precision case,
we suppose that computing~$Ap$ is performed
by running a linearly convergent process, whose rate is given by
$$
\rho=\frac{\sqrt{\lambda_{\max}/\lambda_{\min}} - 1}
          {\sqrt{\lambda_{\max}/\lambda_{\min}} + 1}.
$$ 
This would be the case, for instance, if $A = J^TW^{-1}J$ and only $W$ is known
\footnote{This case occurs in approximately weighted nonlinear least-squares,
for instance in data assimilation for weather forecasting.}. The energy cost of a full accuracy product is then given by $\log(\epsilon_M)/\log(\rho)$,
where $\epsilon_M$ is the IEEE double machine precision,
while that of an inexact product with
accuracy requirement $\omega$ is $\log(\omega)/\log(\rho)$.

For the discrete precision levels case,
we assume that the products $Ap$ can be computed
in double, single, or half precision.
The gain in energy efficiency in this context depend
on the details of the computer architecture.
According to~~\cite{Mats18,High17,GalaHoro11,PuGalaYangSchaHoro16},
a gain between $3\times$ and $5\times$ can be achieved
for each decrease from double to single, and from single to half.
To model this, we assign a unit energy cost for a matrix-vector product in double precision,
a cost $\frac{1}{4}$ for a product in single precision,
and a cost $\frac{1}{16}$ for a product in half precision.

%

In both cases, we sum up the costs over all iterations
to obtain a cost in number of equivalent full double precision
matrix-vector products.

Tables~\ref{res-exact5} to~\ref{res-MM-inexact5-multi}
summarize our results for the accuracy level $\epsilon = 10^{-5}$.
In these tables,

\vspace*{1mm}
\begin{tabular}{lp{5in}}
$\kappa(A)$  & is the condition number of $A$, \\
$n_{it}$     & is the number of iterations required for termination, \\
cost         & is the equivalent number of full accuracy products used, \\
r.res.gap    & is the squared relative residual gap $\frac{1}{2}\normAi{r(x)-r}^2/|q(x_*)|$ at termination, \\
r.sol.err    & is the relative error in the solution value $(q(x)-q(x_*))/|q(x_*)|$ at termination, \\
r.val.err    & is the relative error in the quadratic value $|q(x)-q|/|q(x_*)|$ at termination.
\end{tabular}
\vspace*{1mm}

\noindent
Recall from (\ref{resconds}) and (\ref{quad-decr})
that $\mbox{r.res.gap} \leq \frac{1}{4}\epsilon$
and $\mbox{r.sol.err}\leq \epsilon$,
while r.val.err obeys (\ref{errqval}). 

\subsubsection{Synthetic examples}

In order to illustrate the theory of Section~\ref{rigorous-s},
we first run versions of CG where we use the exact test (\ref{nEj-CG})
(rather than (\ref{app-nEj-CG})), the true $\normAiA{E_j}$,
and the original termination test (\ref{condr}).
This is of course impractical but allows measuring the potential
for inexactness provided by Theorem~\ref{CG-bounds}.
Table~\ref{res-exact5} reports the results obtained
in the continuously-varying accuracy case.
Similar conclusions hold in the discrete precision levels case (not shown).
First, we note that the double precision and inexact variants of CG and CGR
lead to the same relative error on the solution value (r.sol.err).
As targeted, this value is lower than $\epsilon=10^{-5}$
for condition numbers lower than $10^7$ without reorthogonalization. 
Compared to the double precision versions CG and CGR,
the inexact variants iCG and iCGR exhibit very significant potential
savings in the costs of the products $Ap$.
This is especially the case when the conditioning
of the problem is moderate (at most $10^5$).
As it could be argued that the methods discussed
here should be applied on preconditioned systems,
this restriction only moderately affects the
practical applicability of the technique.
Despite an increase in the number of iterations
for iCG compared to CG for condition numbers larger $10^4$,
the overall cost remains lower with iCG.
Finally, the theoretical bound (\ref{nEj-CG}) can be too conservative,
leading to smaller values of r.res.gap than necessary.
  
\begin{table}[htb]

 {\footnotesize
\hspace{-2cm}\begin{tabular}{l|lrrccc|lrrccc}
  method & $\kappa(A)$ & $n_{it}$ & cost & r.res.gap & r.sol.err & r.val.err.
  &$\kappa(A)$ & $n_{it}$ & cost & r.res.gap & r.sol.err & r.val.err. \\
\hline
CG   &       $10^1$  &   11 & 1.1e+01 & 4.6e-30 & 8.3e-07 & 8.8e-16 &       $10^2$ &   34 & 3.4e+01 & 6.0e-29 & 2.0e-06 & 4.6e-16 \\ 
CGR  &        &   11 & 1.1e+01 & 6.4e-30 & 8.3e-07 & 4.4e-16 &        &   34 & 3.4e+01 & 8.6e-29 & 2.0e-06 & 3.2e-15 \\ 
iCG  &        &   11 & 4.0e+00 & 1.1e-08 & 8.5e-07 & 2.9e-07 &        &   34 & 1.3e+01 & 1.0e-09 & 2.0e-06 & 9.1e-07 \\ 
iCGR &        &   11 & 4.0e+00 & 1.1e-08 & 8.4e-07 & 4.8e-06 &        &   34 & 1.3e+01 & 1.0e-09 & 2.0e-06 & 8.2e-07 \\ 
\hline
CG   &      $10^3$   &  104 & 1.0e+02 & 2.6e-27 & 2.4e-06 & 2.1e-07 &    $10^4$     &  313 & 3.1e+02 & 1.4e-25 & 2.4e-06 & 4.0e-08 \\ 
CGR  &        &  104 & 1.0e+02 & 2.7e-27 & 2.3e-06 & 4.1e-15 &        &  263 & 2.6e+02 & 1.3e-25 & 2.5e-06 & 2.5e-14 \\ 
iCG  &        &  105 & 4.3e+01 & 9.8e-11 & 2.4e-06 & 7.7e-07 &        &  323 & 1.5e+02 & 9.9e-12 & 2.4e-06 & 2.1e-07 \\ 
iCGR &        &  104 & 4.3e+01 & 1.1e-10 & 2.3e-06 & 3.1e-07 &        &  263 & 1.2e+02 & 1.2e-11 & 2.5e-06 & 2.1e-07 \\ 
\hline
CG   &    $10^5$    &  928 & 9.3e+02 & 1.0e-23 & 2.5e-06 & 1.1e-07 &   $10^6$     & 2764 & 2.8e+03 & 6.8e-22 & 2.5e-06 & 1.7e-08 \\ 
CGR  &        &  433 & 4.3e+02 & 9.8e-24 & 2.4e-06 & 2.5e-13 &        &  554 & 5.5e+02 & 6.9e-22 & 2.3e-06 & 2.9e-12 \\ 
iCG  &        &  983 & 5.0e+02 & 1.0e-12 & 2.5e-06 & 6.0e-08 &        & 3000 & 1.6e+03 & 3.5e-13 & 2.5e-06 & 9.4e-09 \\ 
iCGR &        &  433 & 2.2e+02 & 3.7e-12 & 2.4e-06 & 6.4e-08 &        &  554 & 3.2e+02 & 5.0e-12 & 2.3e-06 & 2.2e-08 \\ 
\hline
CG   &     $10^7$   & 3000 & 3.0e+03 & 5.4e-20 & 1.3e-02 & 6.0e-07 &  $10^8$      & 3000 & 3.0e+03 & 1.3e-18 & 3.0e-01 & 3.6e-06 \\ 
CGR  &        &  636 & 6.4e+02 & 4.6e-20 & 2.5e-06 & 1.8e-11 &        &  697 & 7.0e+02 & 3.4e-18 & 2.3e-06 & 2.8e-10 \\ 
iCG  &        & 3000 & 1.9e+03 & 1.4e-13 & 1.9e-02 & 5.3e-08 &        & 3000 & 2.1e+03 & 3.2e-14 & 3.4e-01 & 6.2e-06 \\ 
iCGR &        &  636 & 4.0e+02 & 5.7e-12 & 2.5e-06 & 1.6e-07 &        &  697 & 4.7e+02 & 5.6e-12 & 2.3e-06 & 1.0e-07 \\ 
\hline

\end{tabular}
  }

  \caption{\label{res-exact5} Synthetic examples:
	     exact bound in the continuously-varying accuracy case.}
\end{table}

\begin{figure}[htb!]

\begin{center}
\begin{tabular}{cc}
\hspace{-3.5cm}{\footnotesize  Continuously varying precision}     &
\hspace{-1.5cm}{\footnotesize  Discrete precision levels}  \\
\hspace{-3.5cm}\includegraphics[width=4.4in]{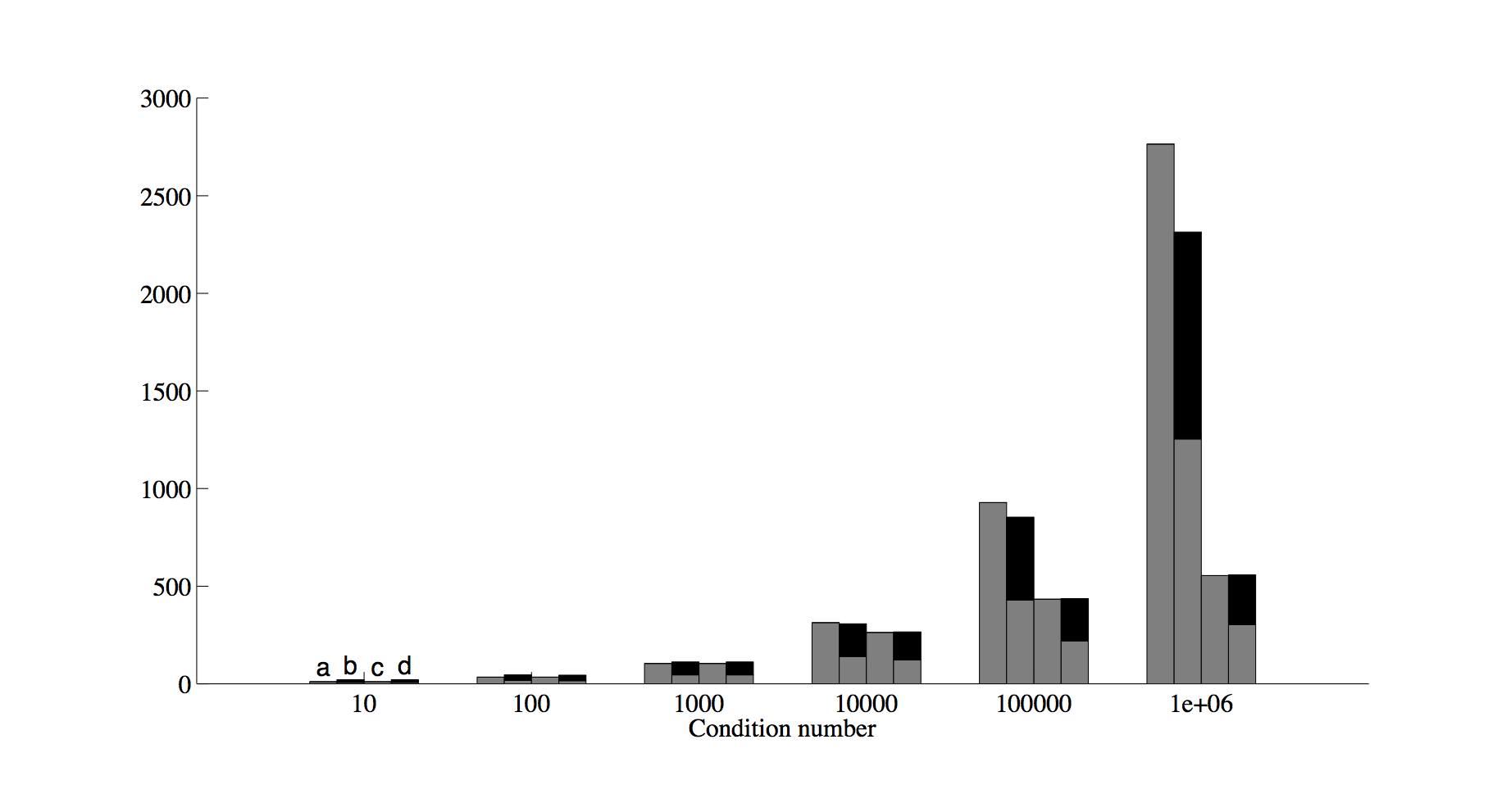} &
\hspace{-1.5cm}\includegraphics[width=4.4in]{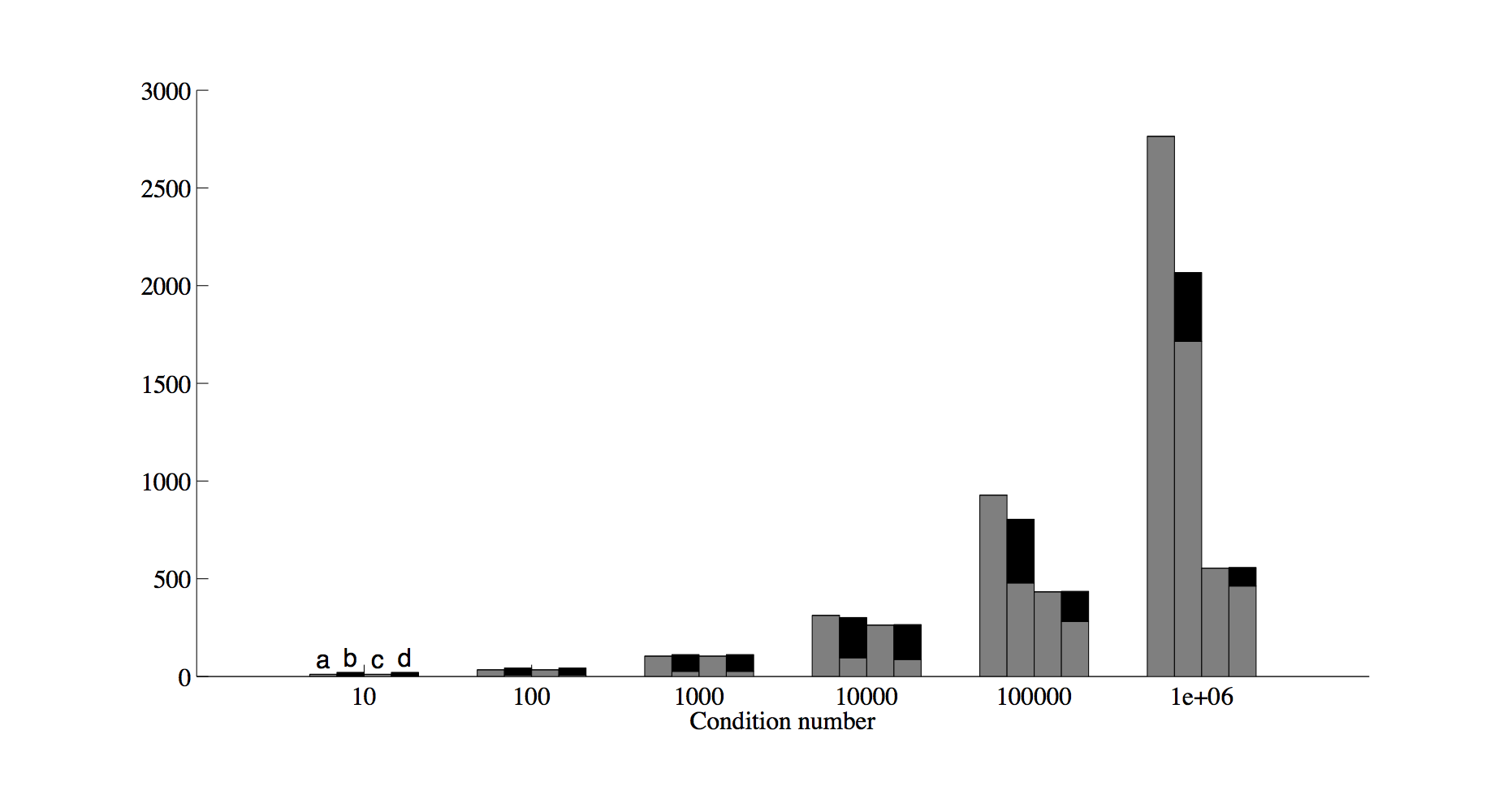} \\
\hspace{-3.5cm}\includegraphics[width=4.4in]{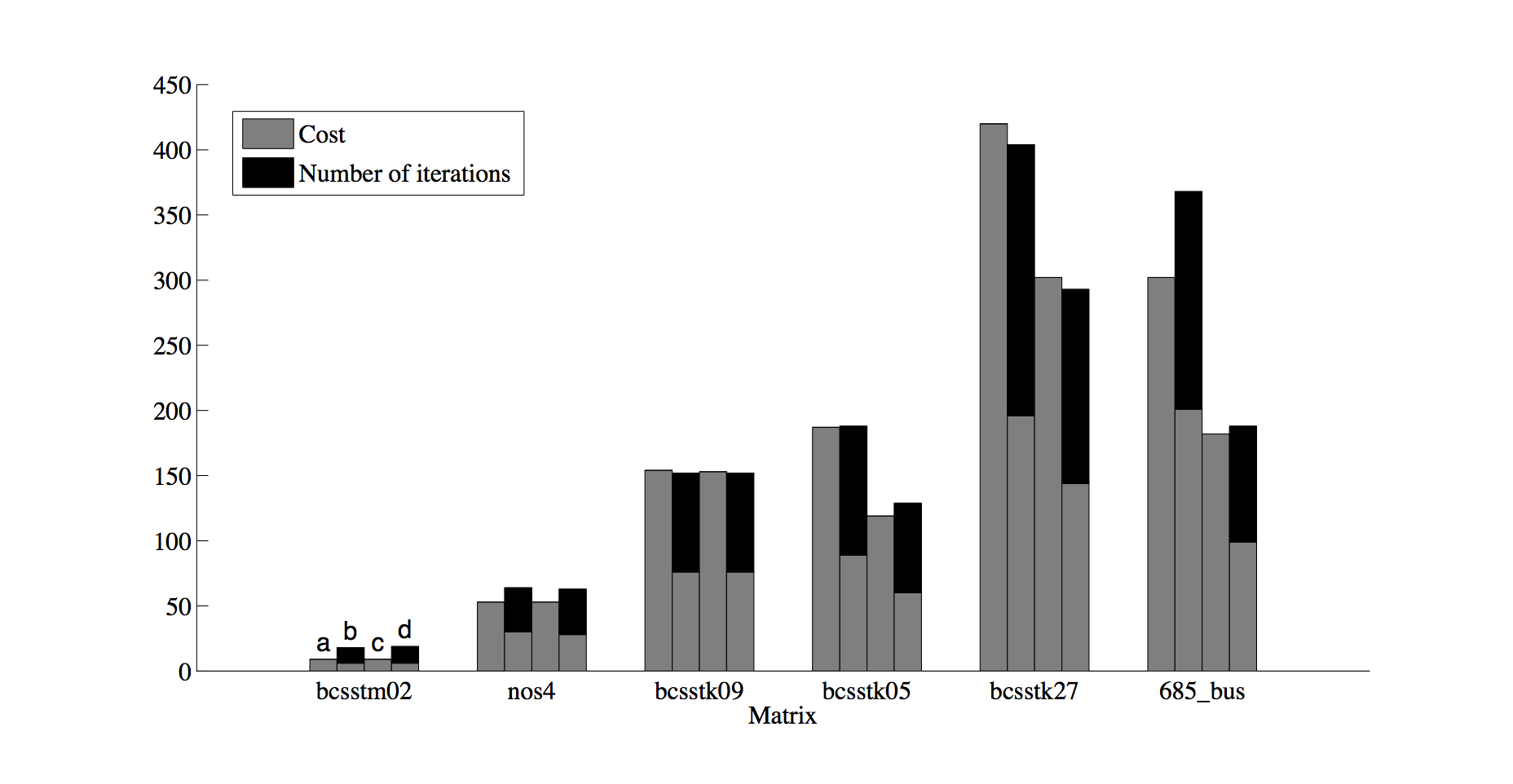} &
\hspace{-1.5cm}\includegraphics[width=4.4in]{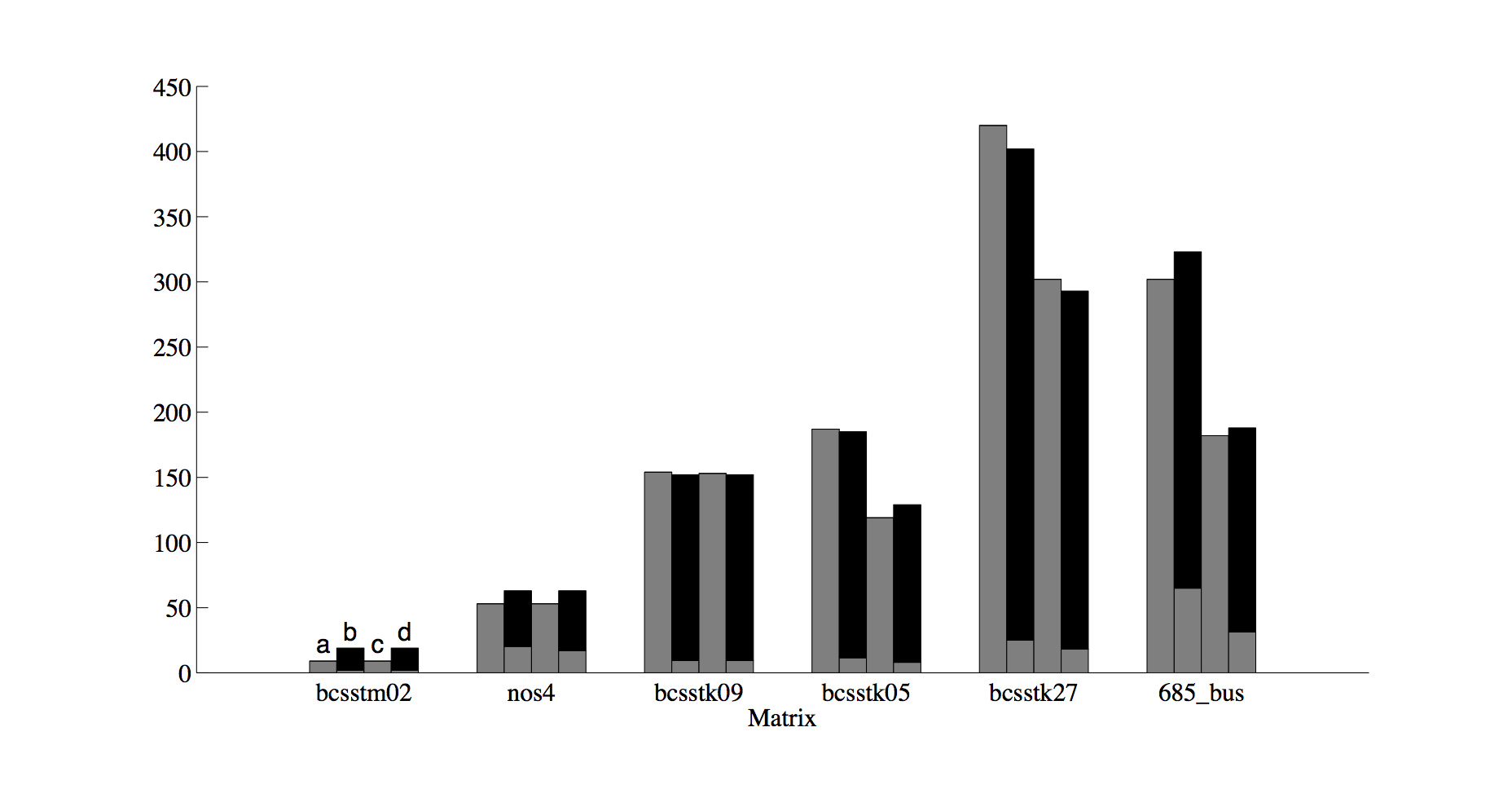} 
\end{tabular}
\caption{Number of iterations and associated costs
for the continuously varying precision and discrete precision levels models.
First row: synthetic matrices with varying logspace condition number.
Second row: matrices from Matrix Market.
(a): double precision CG, (b): inexact practical CG,
(c): double precision CG with reorthogonalization,
(d): inexact practical CG with reorthogonalization.}
\label{fig:ic}
\end{center}

\end{figure}

We now show the effect of using the practical algorithm 
outlined in Section~\ref{practical-CG}.
In addition to using the approximate constants and tests
described in Section~\ref{estimates-s},
we also use estimates of the smallest and largest eigenvalues
obtained by perturbing the true eigenvalues
by a random relative perturbation of magnitude between 0 and 100\%,
with the result that these estimates only hold in order,
but typically have no exact digit.
We report the results of the corresponding runs
in Tables~\ref{res-inexact5-lin} (continuously-varying accuracy case)
and \ref{res-inexact5-multi} (discrete precision levels case),
using the same conventions as for Table~\ref{res-exact5}.
We also provide Figure~\ref{fig:ic} in order to better highlight
the gain in costs obtained using the inexact variants iCG and iCGR for both cases.
In these figures we discarded the cases for which the double precision CG
did not converge before reaching the maximum number of iterations
(high condition numbers).
For each variant, the black bar corresponds to the number of iterations
and the grey bar to the modelled cost.
(For the double precision CG and CGR,
these numbers are equal and only the grey bar is visible.)

\begin{table}[htb]

  {\footnotesize
\hspace{-2cm}\begin{tabular}{l|lrrccc|lrrccc}
  method & $\kappa(A)$ & $n_{it}$ & cost & r.res.gap & r.sol.err & r.val.err.
  &$\kappa(A)$ & $n_{it}$ & cost & r.res.gap & r.sol.err & r.val.err. \\
\hline
CG   &       $10^1$  &   11 & 1.1e+01 & 4.6e-30 & 8.3e-07 & 8.8e-16 &       $10^2$ &   34 & 3.4e+01 & 6.0e-29 & 2.0e-06 & 4.6e-16 \\ 
CGR  &        &   11 & 1.1e+01 & 6.4e-30 & 8.3e-07 & 4.4e-16 &        &   34 & 3.4e+01 & 8.6e-29 & 2.0e-06 & 3.2e-15 \\ 
iCG  &        &   21 & 6.0e+00 & 1.2e-08 & 1.2e-08 & 5.3e-06 &        &   46 & 1.6e+01 & 9.8e-10 & 1.6e-08 & 2.1e-06 \\ 
iCGR &        &   21 & 6.0e+00 & 1.2e-08 & 1.2e-08 & 4.6e-07 &        &   44 & 1.6e+01 & 1.1e-09 & 3.4e-08 & 3.5e-07 \\ 
\hline
CG   &      $10^3$   &  104 & 1.0e+02 & 2.6e-27 & 2.4e-06 & 2.1e-07 &    $10^4$     &  313 & 3.1e+02 & 1.4e-25 & 2.4e-06 & 4.0e-08 \\ 
CGR  &        &  104 & 1.0e+02 & 2.7e-27 & 2.3e-06 & 4.1e-15 &        &  263 & 2.6e+02 & 1.3e-25 & 2.5e-06 & 2.5e-14 \\ 
iCG  &        &  112 & 4.6e+01 & 1.4e-10 & 9.9e-07 & 2.5e-07 &        &  307 & 1.4e+02 & 2.1e-11 & 5.0e-06 & 2.0e-07 \\ 
iCGR &        &  112 & 4.6e+01 & 1.3e-10 & 7.9e-07 & 4.2e-07 &        &  266 & 1.2e+02 & 3.5e-11 & 2.0e-06 & 1.6e-07 \\ 
\hline
CG   &    $10^5$    &  928 & 9.3e+02 & 1.0e-23 & 2.5e-06 & 1.1e-07 &   $10^6$     & 2764 & 2.8e+03 & 6.8e-22 & 2.5e-06 & 1.7e-08 \\ 
CGR  &        &  433 & 4.3e+02 & 9.8e-24 & 2.4e-06 & 2.5e-13 &        &  554 & 5.5e+02 & 6.9e-22 & 2.3e-06 & 2.9e-12 \\ 
iCG  &        &  854 & 4.3e+02 & 4.6e-12 & 1.6e-05 & 7.7e-08 &        & 2314 & 1.3e+03 & 3.4e-12 & 5.6e-05 & 5.0e-07 \\ 
iCGR &        &  436 & 2.2e+02 & 6.5e-11 & 1.9e-06 & 1.9e-09 &        &  558 & 3.0e+02 & 6.5e-10 & 1.6e-06 & 4.2e-07 \\ 
\hline
CG   &     $10^7$   & 3000 & 3.0e+03 & 5.4e-20 & 1.3e-02 & 6.0e-07 &  $10^8$      & 3000 & 3.0e+03 & 1.3e-18 & 3.0e-01 & 3.6e-06 \\ 
CGR  &        &  636 & 6.4e+02 & 4.6e-20 & 2.5e-06 & 1.8e-11 &        &  697 & 7.0e+02 & 3.4e-18 & 2.3e-06 & 2.8e-10 \\ 
iCG  &        & 3000 & 1.8e+03 & 3.4e-12 & 2.0e-02 & 1.6e-06 &        & 3000 & 2.0e+03 & 1.3e-12 & 3.4e-01 & 7.4e-07 \\ 
iCGR &        &  642 & 3.7e+02 & 4.0e-09 & 1.4e-06 & 2.6e-06 &        &  704 & 4.4e+02 & 2.4e-08 & 1.0e-06 & 1.5e-06 \\ 
\hline

\end{tabular}
  }

  \caption{\label{res-inexact5-lin} Synthetic examples:
	    practical algorithms in the continuously-varying accuracy case.}
\end{table}

In the continuous-varying accuracy case,
Table~\ref{res-inexact5-lin} shows that the practical variants of iCG and iCGR
present similar performances to the theoretical inexact CG methods
in terms of errors on the solution value,
except for the practical iCG when the condition number is larger than $10^5$.
For condition numbers in $\{10^5, 10^6\}$,
we note that the number of iterations of the practical iCG is lower
than in double precision CG.
The approximations involved in the stopping criterion (\ref{delay-termination})
lead to early termination, which result in a solution error
slightly larger than the targeted one ($16.10^{-6}$ and $56.10^{-6}$ instead of $10^{-5}$).
For condition numbers larger than $10^7$,
both the double precision CG and inexact iCG
reach the maximum number of iterations without having converged.
Furthermore, the practical methods effectively provide significant gains
in the cost of performing the matrix-vector products in CG,
as highlighted in Table~\ref{res-inexact5-lin} and Figure \ref{fig:ic}.
Finally, we note a slight increase in the number of iterations
for the practical variants of iCG and iCGR at low to medium condition numbers,
which is partly explained by the fact that the termination criterion
is based on the delay $d$ (10 in our case) to assess termination.
The tuning of the parameter $d$ is problem dependent
and should be adapted to the condition number.

Similar conclusions hold for the discrete precision levels case.
Furthermore, the results indicate that the management of the
inaccuracy budget discussed in Section~\ref{inacc_budget-s} is quite
effective. We note that it leads to even more significant efficiency gains for
moderately conditioned problems.
The situation is, however, reversed for the more ill-conditioned cases,
because $\omega$ then exceeds more quickly
the accuracy threshold allowing single precision.
While the small inaccuracy allowed by the bound can be exploited
in the continuous case, this is no longer the case here
and many products are computed in double precision.


\begin{table}[htb]

  {\footnotesize
\hspace{-2cm}\begin{tabular}{l|lrrccc|lrrccc}
  method & $\kappa(A)$ & $n_{it}$ & cost & r.res.gap & r.sol.err & r.val.err.
  &$\kappa(A)$ & $n_{it}$ & cost & r.res.gap & r.sol.err & r.val.err. \\
\hline
CG   &       $10^1$  &   11 & 1.1e+01 & 4.6e-30 & 8.3e-07 & 8.8e-16 &       $10^2$ &   34 & 3.4e+01 & 6.0e-29 & 2.0e-06 & 4.6e-16 \\ 
CGR  &        &   11 & 1.1e+01 & 6.4e-30 & 8.3e-07 & 4.4e-16 &        &   34 & 3.4e+01 & 8.6e-29 & 2.0e-06 & 3.2e-15 \\ 
iCG  &        &   21 & 1.9e+00 & 3.0e-11 & 3.2e-11 & 5.1e-08 &        &   44 & 6.7e+00 & 6.3e-14 & 3.3e-08 & 1.9e-08 \\ 
iCGR &        &   21 & 1.9e+00 & 3.0e-11 & 3.2e-11 & 3.9e-08 &        &   44 & 6.7e+00 & 6.7e-14 & 3.3e-08 & 1.6e-10 \\ 
\hline
CG   &      $10^3$   &  104 & 1.0e+02 & 2.6e-27 & 2.4e-06 & 2.1e-07 &    $10^4$     &  313 & 3.1e+02 & 1.4e-25 & 2.4e-06 & 4.0e-08 \\ 
CGR  &        &  104 & 1.0e+02 & 2.7e-27 & 2.3e-06 & 4.1e-15 &        &  263 & 2.6e+02 & 1.3e-25 & 2.5e-06 & 2.5e-14 \\ 
iCG  &        &  112 & 2.7e+01 & 1.2e-16 & 9.1e-07 & 3.7e-08 &        &  302 & 9.6e+01 & 3.8e-20 & 4.9e-06 & 3.1e-07 \\ 
iCGR &        &  112 & 2.6e+01 & 1.3e-16 & 7.9e-07 & 7.9e-11 &        &  266 & 8.7e+01 & 6.1e-20 & 2.0e-06 & 1.7e-11 \\ 
\hline
CG   &    $10^5$    &  928 & 9.3e+02 & 1.0e-23 & 2.5e-06 & 1.1e-07 &   $10^6$     & 2764 & 2.8e+03 & 6.8e-22 & 2.5e-06 & 1.7e-08 \\ 
CGR  &        &  433 & 4.3e+02 & 9.8e-24 & 2.4e-06 & 2.5e-13 &        &  554 & 5.5e+02 & 6.9e-22 & 2.3e-06 & 2.9e-12 \\ 
iCG  &        &  805 & 4.8e+02 & 2.2e-22 & 2.0e-05 & 3.5e-07 &        & 2067 & 1.7e+03 & 7.2e-22 & 8.9e-05 & 8.1e-07 \\ 
iCGR &        &  436 & 2.8e+02 & 2.1e-21 & 1.9e-06 & 4.9e-14 &        &  558 & 4.6e+02 & 1.1e-21 & 1.6e-06 & 2.0e-12 \\ 
\hline
CG   &     $10^7$   & 3000 & 3.0e+03 & 5.4e-20 & 1.3e-02 & 6.0e-07 &  $10^8$      & 3000 & 3.0e+03 & 1.3e-18 & 3.0e-01 & 3.6e-06 \\ 
CGR  &        &  636 & 6.4e+02 & 4.6e-20 & 2.5e-06 & 1.8e-11 &        &  697 & 7.0e+02 & 3.4e-18 & 2.3e-06 & 2.8e-10 \\ 
iCG  &        & 3000 & 2.9e+03 & 4.4e-20 & 1.3e-02 & 5.9e-07 &        & 3000 & 3.0e+03 & 1.3e-18 & 3.0e-01 & 3.6e-06 \\ 
iCGR &        &  642 & 5.9e+02 & 5.8e-20 & 1.4e-06 & 2.1e-11 &        &  704 & 6.8e+02 & 3.8e-18 & 9.9e-07 & 1.3e-10 \\ 
\hline

\end{tabular}
  }

  \caption{\label{res-inexact5-multi} Synthetic examples:
	     practical algorithms in the discrete precision levels case.}
\end{table}


\subsubsection{Examples from the Matrix Market}

We conclude our experiments with some examples from the NIST Matrix Market. 
Properties of these matrices are given in Table~\ref{mat_prop}. 
All of them are symmetric positive definite and result from discretization of PDEs.
Again, in our computations, we use estimates of the smallest and largest eigenvalues
obtained by perturbating the exact eigenvalues
by a random relative perturbation of magnitude between 0 and 100\%.
We report the results of the corresponding runs
in Tables~\ref{res-MM-inexact5-lin} (continuous-varying accuracy case)
and \ref{res-MM-inexact5-multi} (discrete accuracy levels case),
using the same conventions as for Table~\ref{res-exact5}.
The number of iterations and associated relative costs are also shown in Figure \ref{fig:ic}. 

\begin{table}[htb]
\begin{center}
{
\begin{tabular}{l|lll}
  Matrix  &Dimension& $\kappa_2(A)$ & $\Vert A\Vert_2$  \\
  \hline
  bcsstm02&   66 &   8.8  &   0.17\\
  nos4    &  100 & 1.5e03 &   0.85\\
  bcsstk09& 1083 & 9.5e03 & 6.8e07\\
  bcsstk05&  153 & 1.4e04 & 6.2e06\\
  bcsstk27& 1224 & 2.4e04 & 3.5e06\\
  685\_bus&  685 & 4.2e05 & 2.6e04\\
  nos1    &  237 & 2.0e07 & 2.5e09\\
  nos7    &  729 & 2.4e09 & 9.9e06\\
\hline
\end{tabular}
}
\end{center}
\caption{\label{mat_prop}Properties of the Matrix Market matrices
    (sorted by increasing condition number)}
\end{table}

Similar comments to the synthetic cases can be made.
Regarding the continuous-varying accuracy case,
the practical iCG and iCGR tend to result in an increase of the iteration numbers
compared to both the double precision CG and CGR,
and the impractical methods with inaccurate matrix-vector products
(not shown for the latter ones in Table~\ref{res-MM-inexact5-lin} nor Figure \ref{fig:ic}).
However, a tiny gain in iteration numbers may be obtained
with the practical  iCG and iCGR as observed for the matrix bcsstk27.
Despite these increases in the iteration number,
the costs of the practical methods remain lower
than those of the methods in double precision.
Again the practical algorithms seem more sensitive to rounding errors.
For instance, the desired accuracy is not reached for matrices nos1 and nos7
without reorthogonalization.
As observed in the synthetic matrices,
it can result in a decrease in the number of iterations of
the practical methods, and so a damage of the quality of the solution,
compared to the methods with matrix-vector products in double precision.
Finally, as for the synthetic matrices,
the fact that the delay $d=10$ introduced in the stopping criterion
is too large partly explains the cost increases observed
for matrices with a small condition number (bcsstm02 and nos4).

\begin{table}[htb]

  {\footnotesize
\hspace{-3cm}\begin{tabular}{l|lrrccc|lrrccc}
  method & matrix & $n_{it}$ & cost & r.res.gap & r.sol.err & r.val.err.
  & matrix & $n_{it}$ & cost & r.res.gap & r.sol.err & r.val.err. \\
\hline
CG   &    {\rm bcsstm02}    &    9 & 9.0e+00 & 2.2e-32 & 2.2e-06 & 0.0e+00 &    {\rm nos4}    &   53 & 5.3e+01 & 1.1e-28 & 2.3e-06 & 4.7e-11 \\ 
CGR  &        &    9 & 9.0e+00 & 9.8e-32 & 2.2e-06 & 0.0e+00 &        &   53 & 5.3e+01 & 2.7e-28 & 2.3e-06 & 8.6e-15 \\ 
iCG  &        &   18 & 5.0e+00 & 2.0e-08 & 2.0e-08 & 2.8e-05 &        &   64 & 2.9e+01 & 1.1e-12 & 4.0e-09 & 8.1e-07 \\ 
iCGR &        &    19 & 6.0e+00 & 3.0e-08 & 3.0e-08 & 2.8e-06 &        &   63 & 2.8e+01 & 1.3e-11 & 2.7e-09 & 1.2e-06 \\ 
\hline
CG   &     {\rm bcsstk09}    &  154 & 1.5e+02 & 4.3e-27 & 2.4e-06 & 1.7e-08 &    {\rm bcsstk05}      &  187 & 1.9e+02 & 7.6e-28 & 2.5e-06 & 6.3e-08 \\ 
CGR  &        &  153 & 1.5e+02 & 4.5e-27 & 2.5e-06 & 4.4e-14 &        &  119 & 1.2e+02 & 8.9e-28 & 2.2e-06 & 2.2e-14 \\ 
iCG  &        & 152 & 7.6e+01 & 4.5e-13 & 2.8e-06 & 4.3e-07 &        &  188 & 8.8e+01 & 2.4e-11 & 1.0e-05 & 4.5e-06 \\ 
iCGR &        &  152 & 7.6e+01 & 2.3e-12 & 2.7e-06 & 5.8e-08 &        &  129 & 6.0e+01 & 2.0e-09 & 4.6e-09 & 8.0e-06 \\ 
\hline
CG   &  {\rm bcsstk27}      &  420 & 4.2e+02 & 2.1e-29 & 2.4e-06 & 2.4e-08 &   {\rm 685$\_$bus}     &  302 & 3.0e+02 & 8.2e-27 & 2.5e-06 & 3.3e-08 \\ 
CGR  &        &  302 & 3.0e+02 & 1.9e-27 & 2.4e-06 & 2.9e-15 &        &  182 & 1.8e+02 & 1.4e-26 & 2.3e-06 & 1.4e-14 \\ 
iCG  &        &  404 & 2.0e+02 & 2.6e-12 & 7.7e-06 & 5.0e-08 &        &  368 & 2.0e+02 & 1.6e-13 & 4.8e-06 & 2.3e-07\\ 
iCGR &        &  293 & 1.4e+02 & 4.0e-12 & 4.1e-06 & 7.5e-08 &        &  188 & 9.9e+01 & 3.9e-12 & 1.0e-06 & 4.3e-07 \\ 
\hline
CG   &    {\rm nos1}     &  711 & 7.1e+02 & 3.6e-23 & 3.1e-01 & 6.8e-07 &   {\rm nos7}     & 1810 & 1.8e+03 & 1.3e-19 & 2.1e-06 & 2.6e-08 \\ 
CGR  &        &  220 & 2.2e+02 & 2.5e-23 & 2.1e-06 & 6.1e-13 &        &  270 & 2.7e+02 & 2.9e-19 & 1.8e-06 & 4.0e-12 \\ 
iCG  &        &  711 & 4.4e+02 & 2.2e-12 & 4.0e-01 & 7.3e-06 &        & 1031 & 7.1e+02 & 2.0e-09 & 1.0e-02 & 1.5e-05  \\ 
iCGR &        &  230 & 1.4e+02 & 1.9e-08 & 2.5e-08 & 2.5e-05 &        &  260 & 1.7e+02 & 8.3e-08 & 1.3e-05 & 2.1e-05 \\ 
\hline

\end{tabular}
  }

  \caption{\label{res-MM-inexact5-lin} Matrix Market:
	  practical algorithms in the continuously-varying accuracy case.}
\end{table}

Regarding the discrete precision levels case
(see Table~\ref{res-MM-inexact5-multi} and Figure \ref{fig:ic}),
we note that it leads to even more significant efficiency gains for
moderately conditioned problems and large efficiency gains for the
ill-conditioned cases.  
This can be partly associated with the fact that these
problems are easier than those obtained in the synthetic cases. However, it
remains cases where the desired accuracy cannot be reached without
reorthogonalization due to rounding errors.
As a consequence, a good problem
preconditioning is even more important in the multi-precision context.

\begin{table}[htb]

  {\footnotesize
\hspace{-3cm}\begin{tabular}{l|lrrccc|lrrccc}
  method & $\kappa(A)$ & $n_{it}$ & cost & r.res.gap & r.sol.err & r.val.err.
  &$\kappa(A)$ & $n_{it}$ & cost & r.res.gap & r.sol.err & r.val.err. \\
\hline
CG   &    {\rm bcsstm02}    &    9 & 9.0e+00 & 2.2e-32 & 2.2e-06 & 0.0e+00 &    {\rm nos4}    &   53 & 5.3e+01 & 1.1e-28 & 2.3e-06 & 4.7e-11 \\ 
CGR  &        &    9 & 9.0e+00 & 9.8e-32 & 2.2e-06 & 0.0e+00 &        &   53 & 5.3e+01 & 2.7e-28 & 2.3e-06 & 8.6e-15 \\ 
iCG  &        &   19 & 1.9e+00 & 2.4e-12 & 2.4e-12 & 8.7e-08 &        &   63 & 2.0e+01 & 5.3e-21 & 5.2e-09 & 6.8e-09 \\ 
iCGR &        &   19 & 1.9e+00 & 2.7e-12 & 2.7e-12 & 3.1e-08 &        &   63 & 1.7e+01 & 1.4e-20 & 2.7e-09 & 4.3e-11 \\ 
\hline
CG   &     {\rm bcsstk09}    &  154 & 1.5e+02 & 4.3e-27 & 2.4e-06 & 1.7e-08 &    {\rm bcsstk05}      &  187 & 1.9e+02 & 7.6e-28 & 2.5e-06 & 6.3e-08 \\ 
CGR  &        &  153 & 1.5e+02 & 4.5e-27 & 2.5e-06 & 4.4e-14 &        &  119 & 1.2e+02 & 8.9e-28 & 2.2e-06 & 2.2e-14 \\ 
iCG  &        &  152 & 9.5e+00 & 2.9e-13 & 2.8e-06 & 4.0e-07 &        &  185 & 1.2e+01 & 1.2e-11 & 1.0e-05 & 2.9e-06 \\ 
iCGR &        &  152 & 9.5e+00 & 2.8e-12 & 2.7e-06 & 2.6e-07 &        &  129 & 8.1e+00 & 4.2e-11 & 2.7e-09 & 2.6e-06 \\ 
\hline
CG   &  {\rm bcsstk27}      &  420 & 4.2e+02 & 2.1e-29 & 2.4e-06 & 2.4e-08 &   {\rm 685$\_$bus}     &  302 & 3.0e+02 & 8.2e-27 & 2.5e-06 & 3.3e-08 \\ 
CGR  &        &  302 & 3.0e+02 & 1.9e-27 & 2.4e-06 & 2.9e-15 &        &  182 & 1.8e+02 & 1.4e-26 & 2.3e-06 & 1.4e-14 \\ 
iCG  &        &  402 & 2.5e+01 & 6.5e-13 & 7.7e-06 & 4.3e-08 &        &  323 & 6.5e+01 & 7.5e-18 & 4.7e-06 & 4.6e-08\\ 
iCGR &        &   293 & 1.8e+01 & 6.6e-13 & 4.1e-06 & 7.0e-08 &        &  188 & 3.1e+01 & 7.1e-16 & 1.0e-06 & 3.9e-10  \\ 
\hline
CG   &    {\rm nos1}     &  711 & 7.1e+02 & 3.6e-23 & 3.1e-01 & 6.8e-07 &   {\rm nos7}     & 1810 & 1.8e+03 & 1.3e-19 & 2.1e-06 & 2.6e-08 \\ 
CGR  &        &  220 & 2.2e+02 & 2.5e-23 & 2.1e-06 & 6.1e-13 &        &  270 & 2.7e+02 & 2.9e-19 & 1.8e-06 & 4.0e-12 \\ 
iCG  &        &   711 & 4.5e+01 & 7.0e-15 & 3.9e-01 & 1.7e-06 &        &  957 & 2.1e+02 & 4.4e-15 & 3.0e-03 & 4.9e-08 \\ 
iCGR &        &  230 & 1.4e+01 & 2.5e-12 & 9.0e-09 & 6.1e-07 &        &  269 & 4.1e+01 & 9.6e-14 & 4.1e-06 & 1.3e-08 \\ 
\hline

\end{tabular}
  }

  \caption{\label{res-MM-inexact5-multi} Matrix Market:
	   practical algorithms in the discrete precision levels case.}
\end{table}

\section{Conclusions}\label{concl-s}

We have considered the iterative solution of convex quadratic optimization
problems (\ref{quad-prob}) and linear systems (\ref{system}) 
using the CG algorithm with inaccurate matrix-vector products,
with the aim of monitoring the decrease of the quadratic objective function.
Circumventing the unavailability of some of the quantities involved in the theory,
we have proposed estimates and derived a practical algorithm that use them.
Our numerical experiments suggest that significant gains in energy efficiency
can be achieved by the use of variable precision matrix-vector products. Such gains are most noticeable for problems that are reasonably well-conditioned,
and occur both in the case where the accuracy of the products
can be controlled continuously, and in the case where it is
limited to discrete predefined levels. We have illustrated the latter in the important
context of multi-precision computations. However, the potential speed-up in this context, as well as the gain in energy efficiency, are limited  by the ability of the algorithm to run with precisions lower than the IEEE double precision.  This depends both on the data $(A,b)$ and the decrease in the quadratic that is targeted.
%

In view of the promising potential of this approach,
it may be of interest to apply it in a more general context,
for example, other optimization algorithms involving nonquadratic
and possibly nonconvex objective functions.
It is also worthwhile, in our opinion,
to pursue experimentation with other methods
beyond CG in the framework of multi-precision arithmetic.

While we focused our analysis on inexact matrix-vector products,
a realistic assumption in large-scale applications
where this product often involves the application of
several complicated operators (see \cite{GratGuroSimoToin18b} for example),
the cost of inner-products involved in CG
(and also potentially in reorthogonalization)
may also be significant in some applications.
Strategies to reduce this cost are therefore also of interest.
It is not the purpose of this paper to develop a rigorous analysis
of CG with inexact inner products or reorthogonalization techniques,
but we defer this analysis to a future contribution.

Finally,  it would be interesting to investigate how this approach can be adapted to communication-avoiding algorithms like the $s$-step Krylov methods \cite{vanRos83,ChrGea89}. In recent works on the $s$-step CG, \cite{Car18} suggested a criterion for the adaptive selection of the parameter $s$ (defining the size of the block of iterations) explointing bounds on the residual gap. While combining both approaches would be of interest, this is out of the scope of the present work.

\section*{Acknowledgment}

The authors thank the two anonymous referees for their constructive comments that leads to significant improvements in the presentation. Furthermore, the authors are indebted to Pr. S. Matsuoka (Riken) for an interesting conversation \cite{Mats18} which confirmed their interest in multi-precision arithmetic
in the context of very high performance computing, to Pr A. Podobas (Tokyo Institute of Technology) for providing further pointers
on computer architecture, and to Pr. M. Dayd\'{e} (IRIT) for his continued and friendly support. Ph. Toint was partially supported by ANR-11-LABX-0040-CIMI within the program ANR-11-IDEX-0002-02.

\end{document}